%% file: main.tex
\documentclass[reqno]{siamart220329}
\usepackage[english]{babel}
\usepackage{lipsum}
\makeatletter

\newcommand{\authorfootnotes}{\renewcommand\thefootnote{\@fnsymbol\c@footnote}}%
\makeatother
\usepackage{etoolbox} %

\newcommand*\linenomathpatch[1]{%
  \cspreto{#1}{\linenomath}%
  \cspreto{#1*}{\linenomath}%
  \csappto{end#1}{\endlinenomath}%
  \csappto{end#1*}{\endlinenomath}%
}

\linenomathpatch{equation}
\linenomathpatch{gather}
\linenomathpatch{multline}
\linenomathpatch{align}
\linenomathpatch{alignat}
\linenomathpatch{flalign}

\usepackage{amsfonts}
\usepackage{amssymb}
\usepackage{centernot}
\usepackage{graphicx}
\usepackage{algorithm}
\usepackage{algorithmic}
\usepackage[algo2e]{algorithm2e}
\newtheorem{thm}{Theorem}[section]
\newtheorem{prop}[thm]{Proposition}
\newtheorem{defn}[thm]{Definition}

\usepackage{fullpage}

\usepackage{url}

\title{Learning zeros of Fokker-Planck operators}

\usepackage[displaymath, mathlines]{lineno}

\begin{document}

\maketitle

 {\normalsize
 \centering
  \authorfootnotes
  Pinak Mandal\footnote{\thanks{Corresponding author: \texttt{pinak.mandal@icts.res.in}}}\textsuperscript{1,3}, Amit Apte\textsuperscript{2,3}\par
  \textsuperscript{1} The University of Sydney, NSW 2006 Australia \par
  \textsuperscript{2}Indian Institute of Science Education and Research, Pune 411008 India\par 
  \textsuperscript{3} International Centre for Theoretical Sciences - TIFR, Bangalore 560089 India \par
  \bigskip}

\begin{abstract}
In this paper we devise a deep learning algorithm to find non-trivial zeros of Fokker-Planck operators when the drift is non-solenoidal. We demonstrate the efficacy of our algorithm for problem dimensions ranging from 2 to 10. This method scales linearly with dimension in memory usage. In the problems we studied, overall computational time seems to scale approximately quadratically with dimension. We present results that indicate the potential of this method to produce better approximations compared to Monte Carlo methods, for the same overall sample sizes, even in low dimensions. Unlike the Monte Carlo methods, the deep network method gives a functional form of the solution. We also demonstrate that the associated loss function is strongly correlated with the distance from the true solution, thus providing a strong numerical justification for the algorithm. Moreover, this relation seems to be linear asymptotically for small values of the loss function.
\end{abstract}

\section{Introduction}\label{sec-intro}\input{steady-tex/intro}

\section{Problem statement}
\label{sec-prob}\input{steady-tex/problem}

\section{Previous works}\label{sec-prev-work}\input{steady-tex/prev-work}

\section{Overview of deep learning}
\label{sec-learning}
\input{steady-tex/learning}

\section{The algorithm}
\label{sec-algo}\input{steady-tex/algo}

\section{Results}\label{sec-steady-res}
\input{steady-tex/results}

\section{Conclusions and future work}
\label{sec-conclusions}
\input{steady-tex/conclusions}

\section{Appendix}
\label{sec-appendix}
\input{steady-tex/appendix}

\section*{Acknowledgements}
\input{steady-tex/acknowledge}

\bibliographystyle{siamplain}
\bibliography{ref}
\end{document}

%% file: steady-tex/intro.tex
Many real world problems can be modelled as the response of nonlinear systems to random excitations and such systems have been a topic of interest for a long time. Stochastic differential equations (SDE) provide the natural language for describing many of these systems. Although SDEs have their origins in the study of Brownian motion by Einstein and Smoluchowski, it was Itô who first developed the mathematical theory. Since then SDEs have extensively appeared in physics~\cite{lelievre2016partial, strauss2017hitch, ivanov1980method}, biology~\cite{allen2010introduction}, mathematical finance~\cite{delong2013backward, karoui1997non}, and many other fields~\cite{oksendal2003stochastic, gardiner2009stochastic}. 
The probability density associated with an Itô SDE evolves in time according to a Fokker-Planck equation (FPE) or Kolmogorov forward equation. A stationary FPE (SFPE) can be solved analytically when the corresponding Itô SDE has a drift term that can be represented as the gradient of some potential~\cite{risken1996fokker}. 
But the same is not true when the drift is not of the aforementioned form. In either of these two cases, the time-dependent FPE does not admit a closed form solution in general even when the drift is integrable, thus requiring numerical solutions in most cases. One of the main challenges for numerical solutions of FPE is that the solution of an FPE is a probability density, which requires an integral condition for normalization, which is extremely hard to implement in dimensions larger than two.

In recent times deep learning has been successfully used to solve high-dimensional PDEs~\cite{sirignano2018dgm, han2018solving, raissi2019physics}. Although universal approximation theorems~\cite{pinkus1999approximation, lu2020universal, de2021approximation, kovachki2021universal, de2022error, mishra2022estimates} guarantee existence of neural networks that approximate the true solution well, due to the non-convex nature of loss functions one can not guarantee convergence of neural networks to the true solution during training in many instances~\cite{krishnapriyan2021characterizing, basir2022investigating}.
Moreover, these methods are almost always used for PDEs with boundary conditions not containing integral terms which makes applying them for FPEs challenging. Despite many such issues faced by deep learning solutions to PDEs, it is a worthwhile paradigm to work in while dealing with high-dimensional PDEs for a variety of reasons, some of this are as follows. Most deep learning methods are mesh-free~\cite{blechschmidt2021three} and have the potential to deal with the curse of dimensionality much better than classical methods \cite{cioica2022deep}. Moreover, some of them focus on computing pointwise solutions to PDEs~\cite{han2018solving} which albeit non-standard, might be the only practical and efficient approach in high dimensions. 

The goal of this paper is to devise a reliable, mesh-free deep learning algorithm to solve high-dimensional stationary FPE, i.e., to find non-trivial zeros of Fokker-Planck operators, as explained in detail in the next section. In a sequel~\cite{mandal2024solving}, we devise a method for solving high-dimensional time-dependent FPEs using these zeros. The deep networks we use are the well-known LSTM which are discussed in section~\ref{sec-algo}. We demonstrate the efficacy of this method using examples, described in section~\ref{sec-examples}, where the underlying ODE system possesses a global attractor.
Some of these systems are often used to make simple models in the earth sciences and provide ideal test cases for non-linear filtering algorithms~\cite{carrassi2022data}. The results are discussed in section~\ref{sec-steady-res}. We solve 2, 3, 4, 6, 8, and 10 dimensional problems with our method, in order to explore how our method scales with dimension. We compare our method with Monte Carlo solutions for $d=2$. In the specific examples where an analytical solution is known, we also investigate how the loss function that is minimized during the training of the deep network and the distance of the the true solution from the solution represented by the network are related to each other.

%% file: steady-tex/problem.tex
In this paper we are interested in the stationary 
Fokker-Planck equation
\begin{align}
&\mathcal L p \stackrel{\rm def}{=}
  -\sum_{i=1}^d\frac{\partial(\mu_ip)}{\partial x_i} + \sum_{i=1}^d\sum_{j=1}^d\frac{\partial (D_{ij}p)}{\partial x_i\partial x_j}=0,\quad\mathbf x\in\mathbb R^d\,, \label{eq:SFPE-0}\\
  &\int_{\mathbb R^d}p(\mathbf x)\,d\mathbf x = 1,\quad p(\mathbf x)\ge0\;\;\forall\;\mathbf x\in\mathbb R^d\,. \nonumber
\end{align}
Here $\mu\in C^1(\mathbb R^d; \mathbb R^d)$ is a non-solenoidal vector field, i.e., $\nabla\cdot\mu\not\equiv0$ in $\mathbb R^d$ and $\sigma\in C(\mathbb R^d; \mathbb R^{d\times d})$ is a matrix-valued function such that $D=\frac{1}{2}\sigma\sigma^\top$ is positive-definite.
The operator $\mathcal L$ is known as the \textit{Fokker-Planck operator} (FPO). The goal of this work is to devise an algorithm to find a non-trivial zero of  $\mathcal L$ in a mesh-free manner that works well in dimensions that are challenging for classical PDE-solvers. The motivations behind choosing to find a non-trivial zero of $\mathcal L$ rather than solving~\eqref{eq:SFPE-0}, as well as the motivation for restriction to non-solenoidal vector fields $\mu$, are as follows. 
\begin{itemize}
    \item Numerical integration suffers from the curse of dimensionality~\cite{hinrichs2014curse} and consequently the normalization constraint is extremely challenging to compute in high dimensions.
    \item One of the motivations is to devise an algorithm to solve time-dependent FPEs with unique solutions. We describe in a subsequent work~\cite{dynamicfp2024}, a method to find the normalized solution to the time-dependent FPE that uses a non-trivial zero of $\mathcal L$, even if is may be unnormalized.
    \item When $\mu$ is solenoidal, every constant function is a (unnormalized) zero of $\mathcal L$. In this case, we show in~\cite{dynamicfp2024} that if the corresponding time-dependent FPE has a unique solution, it may be obtained without using a non-trivial zero of $\mathcal L$ to calculate it. Hence we focus on non-solenoidal vector fields in this paper.
    \item Lastly, rather than trying to force normalization during the computation of a non-trivial zero, it is much more economical to integrate the zero at the end to find the normalization constant at a one-time cost. Quasi Monte Carlo~\cite{leobacher2014introduction} or deep learning methods like i-flow~\cite{gao2020flow} can be used for this purpose. 
\end{itemize}

Although our method is valid for any matrix-valued $\sigma$ that gives rise to a positive definite $D$, in the demonstrations we use the form  $\sigma = a I_d$ where $a$ is a positive constant and $I_d$ is the $d\times d$ identity matrix. This allows us to abuse notation and use $\sigma$ and $D$ as scalar quantities. With this simplification our equation becomes, 
\begin{align}
    \mathcal L p= -\nabla\cdot(\mu p) + D\Delta p=0\,.\label{eq:SFPE-1}
\end{align}
where $\Delta$ is the Laplacian operator. 

Since we approximate a solution to~\eqref{eq:SFPE-1} with a neural network, it is sensible to consider strong solutions. We therefore  restrict our search space of functions to $W^{1, 2}_{\rm loc}(\mathbb R^d)\cap C^2(\mathbb R^d)$. Since the superscripts for Sobolev spaces have been used interchangeably in literature, to avoid confusion we define $W^{1, 2}(\mathbb R^d)$ as
\begin{align}
    W^{1, 2}(\mathbb R^d) \stackrel{\rm def}{=} \{f\in L^2(\mathbb R^d):\|\nabla f\|_2\in L^2(\mathbb R^d)\}\,.\label{eq:def-Sobolev}
\end{align} 
Sobolev spaces are frequently encountered while studying elliptic PDEs and therefore are very well-studied~\cite{brezis2011functional, kilpelainen1994weighted}.  This choice of function space enables us to prove uniqueness of solutions to the SFPEs that we will encounter in this paper, as described in detail in appendix~\ref{ssec-unique}. Moreover, density of arbitrary-size neural networks in the space of continuous functions~\cite{pinkus1999approximation} and non-closedness of fixed-size neural networks in Sobolev spaces~\cite{mahan2021nonclosedness} are good justifications for our algorithm, as discussed in greater detail in section~\ref{ssec-infinite-to-finite}, making $W^{1, 2}_{\rm loc}(\mathbb R^d)\cap C^2(\mathbb R^d)$ an ideal function space to work with.

As mentioned earlier, one of the motivations for studying the SFPE is to device methods to find solutions of time-dependent FPE which occur quite often when studying time evolution of probability densities of random dynamical systems governed by stochastic differential equations of the form
\begin{align}
    x(t) = x(0) + \int_0^t \mu(x(s)) ds + \int_0^t \sigma(x(s)) dW(s) \,,
\label{eq-sde}\end{align} 
where $W(s)$ is a $d$-dimensional Brownian motion. We also refer to the related ODE system $\dot{x} = \mu(x)$ as the associated deterministic system and discuss the relation between the properties of this ODE, in particular, presence of a global attractor, and the solution of the SPFE.

\section{Examples}\label{sec-examples}
From an algorithmic perspective, it is important to have access to a class of equations on which our algorithm can be validated easily. Since classical methods do not work satisfactorily for our problem dimensions, the validating examples we use are those for which the analytical solutions are known. In addition to these examples, we also present results for which analytical solutions are not known. These example problems are described in this section.

\subsection{Gradient systems}
We first describe a very large class of systems for which analytical solutions can be written down, so that they can be used to validate the algorithms we propose. This class consists of equations where the drift $\mu$ can be written as the gradient of a potential function,
\begin{equation}
    \mu = -\nabla V \label{eq:mu-grad} \,.  
\end{equation}
It is easy to verify that $p$ given below is a solution in this special case.
\begin{align}
   p = c\exp\left(-\frac{V}{D}\right) \,.\label{eq:grad-sol}
\end{align}
We refer to a system satisfying~\eqref{eq:mu-grad} as a \textit{gradient system}.
In this paper we use the following gradient systems to validate our algorithm in high dimensions. 

\subsubsection{2D ring system}
For $d=2$, $V=(x^2+y^2-1)^2$, and $\mu=-\nabla V$, we get the following SFPE,
\begin{align}
4(x^2+y^2-1)\left(x\frac{\partial p}{\partial x}+y\frac{\partial p}{\partial y}\right) + 8(2x^2+2y^2-1)p + D\Delta p=0 \,. \label{eq:ring2D}
\end{align}
This system possesses a unique solution concentrated around the unit circle. The proof of uniqueness using the method of Lyapunov functions is given in the appendix~\ref{sssec-2D-unique}. The corresponding ODE system has the unit circle as a global attractor. This is a recurring theme in all of our example problems. Such systems with attractors are of great interest in the study of dynamical systems~\cite{ott1981strange} as well as filtering theory~\cite{kontorovich2009non}. We solve this system for $D=1$.

\subsubsection{2nD ring system}\label{ssec-2nD-ring} We can daisy-chain the previous system to build decoupled systems in higher dimensions. In this case the potential is given by
\begin{align}
    V(\mathbf x)=\sum_{i=0}^{\frac{d}{2}-1}(x_{2i}^2+x_{2i+1}^2-1)^2 \,,\qquad d=2n \,.\label{eq:2nD-V}
\end{align}
Since our algorithm does not differentiate between coupled and decoupled systems, this example serves as a great high-dimensional test case. In a subsequent work~\cite{dynamicfp2024}, we show how to solve the time-dependent FPEs with a method that is intimately related to the method presented in this paper, and this system, being a decoupled high-dimensional system, presents a great way to verify the time-dependent algorithm. This is important since analytical solutions for time-dependent FPEs are not known in general even for gradient systems. Uniqueness of solution for the 2nD ring system directly follows from the uniqueness of solution for the 2D ring system, again thanks to its decoupled nature. In this paper, we solve this system for $d = 2, 4, 6, 8, 10$ and with the choice of diffusion parameter $D=1$.

\subsection{Non-gradient Systems}
Not all the drifts $\mu$ can however be represented as the gradient of a potential. We call the systems belonging to this complementary class, \textit{non-gradient systems}. Analytic solutions for these systems are not known in general. Two examples of such systems that we use in this paper are described below.

Since analytic solutions for non-gradient systems are not known, we restrict our attention to $d=3$ in this case. This is a dimension that can be reliably tackled with Monte Carlo simulations for comparison at a reasonable computational cost. See~\ref{ssec-MC-algo} for a description of the Monte Carlo algorithm that we use for such comparisons.

\subsubsection{Noisy Lorenz-63 system}
\begin{figure}[t!]
    \centering\includegraphics[scale=0.55]{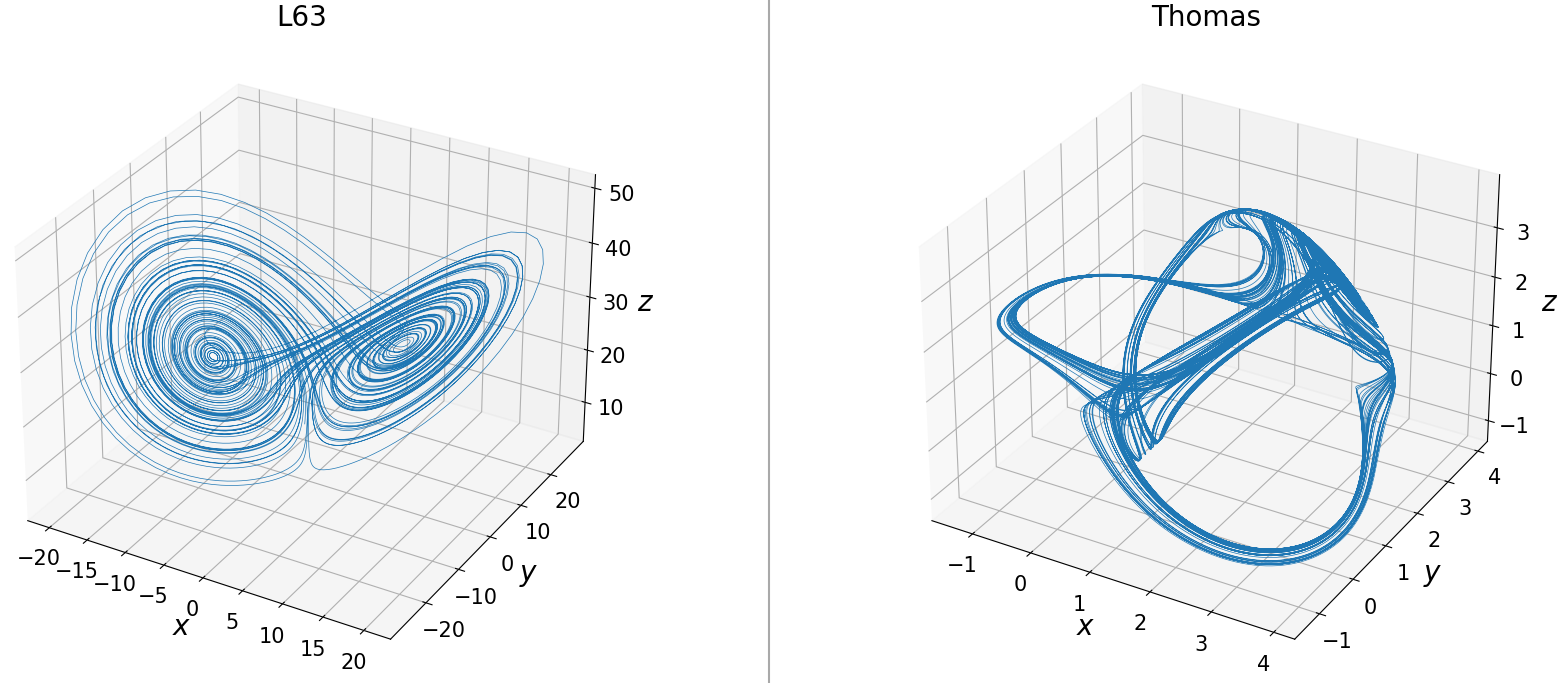}
   \caption{Attractors for non-gradient examples} \label{fig:attractors}
\end{figure}

A famous example is the Lorenz-63 ODE system, first proposed by Edward Lorenz~\cite{lorenz1963deterministic} as an oversimplified model for atmospheric convection, with the drift $\mu$ given by
\begin{align}
    &\mu=[\alpha (y-x),\, x(\rho-z) - y,\, xy - \beta z]^\top \,, \label{eq:mu-L63}\\
    &\alpha = 10 \,, \, \beta = \frac{8}{3}\,, \, \rho=28 \,.
\end{align}
The well-known butterfly attractor associated with the corresponding ODE is shown in figure~\ref{fig:attractors}. 
This system and its variants like Lorenz-96 have become staple test problems in the field of data assimilation~\cite{carrassi2022data, yeong2020particle}. We use the standard parameters to define the drift and solve the system for $D=50$. This choice is motivated by the fact that this system already appears as a test case in~\cite{chen2018efficient}. With the choice of the drift $\mu(x)$ given above, the SFPE~\eqref{eq:SFPE-0} has a unique solution -- for a proof see appendix~\ref{sssec-L63-unique}.

\subsubsection{Noisy Thomas system}
Another example of a non-gradient system that we study is one for which the deterministic version was proposed by René Thomas~\cite{thomas1999deterministic}. It is a 3-dimensional system with cyclical symmetry in the three coordinates $x, y, z$ and the corresponding ODE system has a strange attractor which is depicted in figure~\ref{fig:attractors}.
We solve this system for $D=1$. The SFPE for this problem also has a unique solution -- for a proof see appendix~\ref{sssec-Thomas-unique}.
\begin{align}
    &\mu=[\sin y - bx,\, \sin z - by,\, \sin x - by ]^\top \,, \qquad b =  0.2 \,. \label{eq:mu-Thomas} 
\end{align}

%% file: steady-tex/prev-work.tex
An extensive amount of work has been done on the topic of numerically solving Fokker-Planck equations. A large numer of these works are based on traditional PDE solving techniques such as finite difference~\cite{berezin1987conservative, whitney1970finite, sepehrian2015numerical} and finite element~\cite{naprstek2014finite, masud2005application} methods. For a comparison of these traditional methods we refer the reader to the comparative study~\cite{pichler2013numerical} by Pitcher et al.~where the methods have been applied to 2 and 3 dimensional examples.

In recent times efforts have been made to devise methods that are applicable in dimensions higher than 3. Tensor decomposition methods~\cite{Hackbusch2005HierarchicalKT, kolda2009tensor} are an important toolkit while dealing with high-dimensional problems and they are proving to be useful in designing numerical solvers for PDEs~\cite{ballani2013projection, kressner2010krylov}.  For stationary Fokker-Planck equations, Sun and Kumar~\cite{sun2014numerical} proposed a tensor decomposition and Chebyshev spectral differentiation based method. In this method drift functions are approximated with a sum of functions that are separable in spatial variables, an well-established paradigm for solving PDEs. The differential operator for the stationary FPE is then discretized and finally a least squares problem is solved to find the final solution. The normalization is enforced via addition of a penalty term in the optimization problem. The high-dimensional integral for the normalization constraint in this method is replaced with  products of one dimensional integrals and therefore becomes computable.   

Chen and Majda~\cite{chen2018efficient} proposed a hybrid method that utilizes both kernel and sample based density approximation to solve FPEs that originate from a specific type of SDE referred to as a conditional Gaussian model.
The special structure of the SDE allows one to approximate 
the marginal of a subset of variables as a Gaussian mixture with parameters that satisfy auxiliary SDEs while the marginal of the remaining variables 
is approximated with a non-parametric kernel based method. Finally the joint distribution 
is computed with a hybrid expression. Using this method Chen and Majda computed the solution to a 6 dimensional conceptual model for turbulence. Note that, among our examples only L63 falls under this special  structure.

In recent years machine learning has also been applied to solve SFPEs. Xu et al~\cite{xu2020solving} solved two and three dimensional stationary FPEs with deep learning. Their method enforced normalization via a penalty term in the loss function that represented a Monte-Carlo estimate of the solution. Although simple and effective in lower dimensions, this normalization strategy loses effectiveness in higher dimensions. Zhai et al~\cite{zhai2022deep} have proposed a combination of deep learning and Monte-Carlo method to solve stationary FPEs. The normalization constraint here is replaced with a regularizing term in the loss function which tries to make sure the final solution is close to a pre-computed Monte-Carlo solution. This strategy is more effective than having to approximate high-dimensional integrals and the authors successfully apply their method on Chen and Majda's 6 dimensional example.

%% file: steady-tex/learning.tex
In this section we describe the general process of \textit{learning} a solution to a partial differential equation. The strategy described here will be an integral part of the final algorithm. Machine learning solutions to PDEs can refer to any of the many different scenarios such as super-resolution or using classical grid-based solutions to approximate solutions on finer grids~\cite{li2020fourier}, extension of domain or using classical numerical solutions to approximate solutions on previously unexplored domains~\cite{ovadia2023ditto}, interpolation or extrapolation in parameters of the PDEs i.e. approximating solutions for previously unexplored parameters of the PDE~\cite{li2020neural, khoo2021solving}, learning pointwise solutions to PDEs using associated SDEs \cite{han2018solving}, learning solutions to PDEs globally with a functional form \cite{sirignano2018dgm} etc. Note that, in a lot of these scenarios one uses supervised learning to extend solutions computed with classical methods. These setups are therefore suitable only for low-dimensional problems. Since one of our goals is to solve high-dimensional problems for which classical numerical methods are prohibitively time or memory-consuming, in this work we use unsupervised learning where we do not have access to any pre-computed solutions.

Most previous 
works~\cite{blechschmidt2021three, sirignano2018dgm, yu2018deep, raissi2019physics} in this unsupervised scenario deal with PDEs with boundary conditions, of the type~\eqref{eq:generic-pde}, and hence this section is focused on such equations, rather than our problem~\eqref{eq:SFPE-0} which has no boundary conditions but a normalization condition. 
In particular, we discuss in this section the general \textit{physics-informed} methodology that has been used 
to solve a generic time-independent PDE~\eqref{eq:generic-pde} with a Dirichlet boundary 
condition.
The interested reader can see~\cite{raissi2019physics, blechschmidt2021three, sirignano2018dgm} for more discussions. In the next few subsections we keep simplifying our PDE problem until it finally becomes solvable on a computer.

\subsection{From PDE to optimization problem} In the context of machine learning, \textit{learning} refers to solving an optimization problem. So to solve our PDE with deep learning we first transform it into an optimization problem. For this purpose, we recall the 2nd order PDE we want to solve can be written as 
\begin{equation}
\begin{aligned}
    &\mathcal Lf(\mathbf x) = 0,\quad \mathbf x\in\Omega \,, \qquad \textrm{and} \qquad
    &f(\mathbf x) = g(\mathbf x),\quad \mathbf x\in\partial\Omega \,,
\end{aligned}\label{eq:generic-pde}
\end{equation}
and we are interested in finding a solution in $W^{1,2}_{\rm loc}(\Omega)\cap C^2(\Omega)$. Instead of trying to solve \eqref{eq:generic-pde} a popular strategy is to try to solve the following problem (see for example \cite{sirignano2018dgm}),
\begin{align}
    f^* = \underset{f\in W^{1,2}_{\rm loc}(\Omega)\cap C^2(\Omega)}{\rm arg\,inf}\left[\int_\Omega (\mathcal L f)^2 + \int_{\partial\Omega}(f-g)^2\right] \,.\label{eq:generic-pde-opt}
\end{align}
The choice of function space ensures one-to-one correspondence between the solutions of the PDE and the optimization problem.

\subsection{From infinite-dimensional search space to finite-dimensional search space}\label{ssec-infinite-to-finite} To solve a problem on a machine with finite resources we need to find a finite dimensional (and in fact, a finite) approximation of the infinite dimensional aspects of the problem. We then solve the approximate, finite problem and preferably also estimate how well the finite solution approximates the solution to the original problem.

In particular, we will replace our search space $W^{1,2}_{\rm loc}(\Omega)\cap  C^2(\Omega)$ with a finite dimensional one by appealing to universal approximation theorems that say that neural networks of even the simplest architectures are dense in continuous functions, see for example theorem 3.2 in \cite{kidger2020universal} or proposition 3.7 in \cite{pinkus1999approximation}. Universal approximation theorems typically allow networks to have either arbitrary depth or arbitrary width in order to achieve density \cite{pinkus1999approximation}, \cite{de2021approximation}. But the sets of neural networks with arbitrary depth or width are still infinite dimensional and therefore are infeasible to work with. In practice, we fix an architecture $\mathcal A$ with a fixed number of layers and trainable parameters and work with the following set instead.
\begin{align}
    S_{\mathcal A}\stackrel{\rm def}{=}\{n^{\mathcal A}_\theta: \theta\in \mathbb R^C\}\,.\label{eq:search-space-net}
\end{align}
Here $n^{\mathcal A}_\theta$ is a network with architecture $\mathcal A$ with trainable parameters $\theta$ and $C$ is the total number of trainable parameters or the size of $\theta$. Since $C$ is fixed, $S_{\mathcal A}$ {has a one-to-one correspondence with} $\mathbb R^C$ and therefore is finite-dimensional. Note that, one can impose restrictions on the parameters of the network for example, by regularizing them \cite{shen2022consistency} or fixing a subset of them \cite{zhuang2020comprehensive}. In such cases $S_{\mathcal{A}}$ might only have a one-to-one correspondence with a subset of $\mathbb R^d$ but we do not consider such cases in this work.  Even though we lose the density argument while working with $\theta$ of fixed size, in recent times it has been shown that sets like $S_\mathcal A$ are not closed in $W^{1, 2}(\Omega)$ and $n^{\mathcal A}_\theta$ can be used as a good function approximator, see section~3 in \cite{mahan2021nonclosedness} for a detailed discussion. In the following discussion we suppress the architecture and use $n^{\mathcal A}_\theta$ and $n_\theta$ interchangeably for notational convenience. After restricting our search space to \eqref{eq:search-space-net},  our optimization problem becomes
\begin{align}
    \theta^* = \underset{\theta\in\mathbb R^C}{\rm arg\,inf}\left[\int_{\Omega} (\mathcal L n_\theta)^2 + \int_{\partial\Omega}(n_\theta-g)^2\right] \label{eq:generic-pde-opt-theta}
\end{align}
and the corresponding $n_{\theta^*}$ approximates the solution $f^*$ to the problem~\eqref{eq:generic-pde-opt}.

\subsection{From integrals to sums} 
When the domain $\Omega$ is high dimensional, computation of the integrals in \eqref{eq:generic-pde-opt-theta} will be extremely challenging. To deal with this we will replace the integrals in \eqref{eq:generic-pde-opt-theta} with Monte-Carlo sums.
\begin{align}
    \theta^* = \underset{\theta\in\mathbb R^C}{\rm arg\,inf}\left[\frac{1}{N}\sum_{j=1}^N (\mathcal L n_\theta(\mathbf x_j))^2 + \frac{1}{M}\sum_{j=1}^M(n_\theta(\mathbf y_j)-g(\mathbf y_j))^2\right]\,,\label{eq:generic-opt-final}
\end{align}
where $\{\mathbf x_j\}_{j=1}^N$, $\{\mathbf y_j\}_{j=1}^M$ are uniform samples from $\Omega$ and $\partial \Omega$ respectively. In this case, \eqref{eq:generic-opt-final} can interpreted as trying to find a network that satisfies the original problem~\eqref{eq:generic-pde} at the specified points $\{\mathbf x_j\}_{j=1}^N, \{\mathbf y_j\}_{j=1}^M$, which we can refer to as \textit{collocation points}. 

\subsection{Finding the optimal parameters}\label{ssec-finding-theta}
Having transformed the problem~\eqref{eq:generic-pde-opt} to the one stated above, we perform gradient descent with respect to $\theta$ to find the optimal network for the problem \eqref{eq:generic-opt-final}. The Monte-Carlo sample sizes is dictated by the hardware available.  In our experiments $N = M = 1000$. In many cases, these choices may not be enough to approximate the original integrals sufficiently well, as is the case in the examples in this paper and in general in most problems of interest. In order to overcome this limitation and in order to the learn the solution on the entire domain as thoroughly as possible, we resample the domain every few training iterations. Thus, even though we are limited in sample size by our hardware, we can shift the burden on space or memory to time or number of training iterations, in order to adequately sample the entire domain. This principle of space-time trade-off is ubiquitous in machine learning \cite{buduma2022fundamentals} and comes in many different flavours like mini-batch gradient descent, stochastic gradient descent etc. Even though in this paradigm we are not training our network with typical input-output pairs, our method can be thought of as a variant of the mini-batch gradient descent.

\subsection{Rationale for deep learning}\label{ssec-rationale}
In this context of our problem, deep learning refers to learning an approximate solution to \eqref{eq:generic-pde} with the outlined method with an architecture $\mathcal A$ that is \textit{deep} or has many hidden layers. Deep networks are more efficient as approximators than shallow networks in the sense that they require far fewer number of trainable parameters to achieve the same level of approximation. For a discussion see section 5 of \cite{holstermann2023expressive} or section 4 of \cite{lu2017expressive}. Now that we have described the general procedure of \textit{deep learning} a solution to a PDE, we will pause briefly to point out some benefits and demerits of this approach. Deep learning has, like any other method some disadvantages. 
\begin{itemize}
    \item Deep learning is slower and less accurate for lower dimensional problems for which standard solvers exist and have been in consistent development for many decades. 
    \item Most modern GPUs are optimized for computation with single precision or float32 numbers and float64 computations on GPU are significantly slower than float32. Lower precision float32 is efficient and sufficient for rendering polygons or other image processing tasks which are the primary reasons GPUs were invented~\cite{peddie2023history} but float32 might not be accurate enough for scientific computing.
    \item The objective or \textit{loss} function used in a typical problem might not be convex and hence difficult to deal with~\cite{krishnapriyan2021characterizing, basir2022investigating}, due to multiple local minima. 
\end{itemize}
But even with these disadvantages, the benefits of deep learning make it a worthwhile tool for solving PDEs.
\begin{itemize}
    \item Since we don't need to deal with meshes or grids in this method, we can mitigate the curse of dimensionality in memory. It will be clear from our experiments that the size of the network $C$ does not need to grow exponentially with the dimensions. This method lets one compute the solution at collocation points but if one wants to compute the solution over the entire domain, one needs to sample the entire domain thoroughly which can be done in a sequential manner without requiring more memory as discussed in~\ref{ssec-finding-theta}.
    \item All derivatives are computed with automatic differentiation and therefore are accurate up to floating point errors. Moreover, finite difference schemes do not satisfy some fundamental properties of differentiation e.g. the product rule \cite{ranocha2019mimetic}. With automatic differentiation one does not have to deal with such problems. 
    \item If one computes the solution over the entire domain, the solution is obtained in a functional form which can be differentiated, integrated etc.
    \item Other than a method for sampling no modifications are required for accommodating different domains. 
\end{itemize}

%% file: steady-tex/algo.tex
In this section we outline the algorithm for learning zeros of FPOs. But before that we go through the primary challenges and ways to mitigate them.

\subsection{Unboundedness of the problem domain}\label{ssec-unbounded-domain} We can try the same procedure as outlined in section~\ref{sec-learning} to find a non-trivial zero of $\mathcal L$. But computationally we can only deal with a bounded domain. Hence we focus on a compact domain which contains most of the mass of the solution to \eqref{eq:SFPE-0}. We refer to this domain as the \textit{domain of interest} $\Omega_I$ in the following discussion. {We note that the support of non-trivial zeros of $\mathcal L$ will usually be unbounded and of course we do not know the domain that may contain most of the mass. Thus the choice of $\Omega_I$ needs to be informed by some a priori knowledge about the solution, which in the examples we discuss is related to some attracting set of the deterministic system associated to the drift term $\mu$, i.e., the first term in~\eqref{eq-sde}. We do not need a precise knowledge of such an attracting set. But the smaller the domain $\Omega_I$, the more efficient the proposed method will be, which requires uniform samples from $\Omega_I$.}

\subsection{Existence of the trivial solution}\label{ssec-exist-0} {Since $\mathcal L$ is a linear operator, zero is a trivial solution: $\mathcal L 0 = 0$. We also note that if $\nabla \cdot \mu \not\equiv 0$, then no other constant function is a zero of $\mathcal L$.} Since we want to find a non-trivial zero of $\mathcal L$, we would like avoid the learning the zero function during the training of the network. To deal with this problem \cite{zhai2022deep} added a regularization term that used approximate solutions of \eqref{eq:SFPE-0} found using Monte-Carlo. Here we propose a method that does not require a priori knowing an approximate solution. Consider the operator $\mathcal L_{\rm log}$ instead.
\begin{align}
    \mathcal L_{\rm log}f \stackrel{\rm def}{=} e^{-f}\mathcal L e^f\,.\label{eq:def-log-FPO}
\end{align}
Note that if $f$ is a zero of $\mathcal L_{\rm log}$, then $p = e^f$ is a zero of $\mathcal{L}$, which automatically assures positivity of the solution. Thus we can look for a zero of $\mathcal L_{\rm log}$ to find a non-trivial zero of $\mathcal L$. Straightforward calculation yields
\begin{align}
    \mathcal L_{\rm log}f=
    -\nabla\cdot \mu - \mu \cdot \nabla f + D\left(\|\nabla f\|_2^2 + \Delta f\right)\,.\label{eq:log-FPO}
\end{align}
We again note that when $\nabla\cdot\mu\not\equiv0$, then any constant function can not be a zero of $\mathcal L_{\rm log}$.

\subsection{The steady state algorithm}
The procedure outlined in section~\ref{sec-learning} together with the modifications in sections~\ref{ssec-unbounded-domain}-\ref{ssec-exist-0} immediately yield the following loss function. 
\begin{align}
    L_{\log}(\theta) = \frac{1}{N}\sum_{i=1}^N\mathcal L_{\rm log}(n_\theta(\mathbf x_i))^2\,,\label{eq:def-steady-loss}
\end{align}
where $\{\mathbf x_i\}_{i=1}^N$ is a uniform sample from $\Omega_I$. Accordingly, the final procedure for finding a non-trivial zero of $\mathcal L$ is given in algorithm~\ref{algo:steady}.
\begin{algorithm}[t!]
Select the desired architecture for $n_\theta$.\\
Select resampling interval $\tau$.\\
Select an adaptive learning rate $\delta(k)$ and the number of training iterations $E$. 
Sample $\{\mathbf x_i\}_{i=1}^N$ from $\Omega_I$, the domain of interest.\\
\For {$k=1,2\cdots, E$}{
Compute $\nabla_\theta L_{\log}=\frac{1}{N}\sum_{i=1}^N\nabla_\theta(\mathcal{L}_{\log}(n_\theta(\mathbf x_i))^2)$\\
where $\mathcal L_{\log}f = -\nabla\cdot \mu - \mu \cdot \nabla f + D\left(\|\nabla f\|_2^2 + \Delta f\right)$\\
Update $\theta\leftarrow\theta - \delta(k) \nabla_{\theta}L_{\log}$\\
\If{$k\text{ is divisible by }\tau$}{Resample $\{\mathbf x_i\}_{i=1}^N$ from $\Omega_I$}
}
$e^{n_\theta(\mathbf x)}$ is a non-trivial zero of $\mathcal{L}$.\\
Optional: Approximate $Z\leftarrow\int_{\mathbb R^d}e^{n_\theta(\mathbf x)}\,d\mathbf x$.\\
$\frac{1}{Z}e^{n_\theta(\mathbf x)}$ is the learned, normalized steady state.
\caption{The steady state algorithm}
\label{algo:steady}
\end{algorithm}



We note that such iterative algorithms can use multiple stopping criteria, in addition to number of iterations. Three most common ones use pre-chosen thresholds for the following quantities: (i) $\nabla_\theta L_{\log}$, (ii) the change in $L_{\log}$ with respect to iterations, or (iii) the loss $L_{\log}$ itself. 
Out of these three, the third one is an ineffective criterion in our problem since the threshold value for stopping the algorithm will depend on the network size and architecture. Thus it is difficult to choose the threshold for $L_{\log}$ a priori. Additionally, these criteria can be effective only when the appropriate nuances in their implementation are considered. For example, the changes in the loss, as a functional, every $\tau$ iterations due to domain resampling may need to be taken into account. A more extensive study of the utility of these stopping criteria and the choice of associated hyperparameters would be an interesting direction for further investigations, but in this paper we report the results with the algorithm being run for a pre-chosen number of iterations.

In the following sections we describe in detail the network architecture and optimizer used in our experiments. 

\subsection{Architecture}\label{ssec-architecture}
We choose the widely used LSTM~\cite{sherstinsky2020fundamentals, vennerod2021long} architecture described below for our experiments. This type of architecture rose to prominence in deep learning because of their ability to deal with the vanishing gradient problem, see section IV of \cite{sherstinsky2020fundamentals}, section 2.2 of \cite{vennerod2021long}. A variant of this architecture has also been used to solve PDEs \cite{sirignano2018dgm}. This kind of architectures have been shown to be universal approximators \cite{schafer2006recurrent}. We choose this architecture simply because of how \textit{expressive} they are. By expressivity of an architecture we imply its ability to approximate a wide range of functions and experts have attempted to formalize this notion in different ways in recent times~\cite{lu2017expressive,raghu2016survey,  raghu2017expressive}. Some architectures are probability densities by design i.e. the normalization constraint in~\eqref{eq:SFPE-0} is automatically satisfied for them, see for example~\cite{uria2013rnade, papamakarios2019neural}. But our experiments suggest these architectures are not expressive enough to learn solutions to PDEs efficiently since the normalization constraint makes their structure too rigid. These are the main reasons we choose to focus on learning a non-trivial zero of $\mathcal L$ rather than solving \eqref{eq:SFPE-0}, using LSTM networks which are expressive enough to solve all the problems listed in section~\ref{sec-examples}. We note that such networks have been used in a wide variety of problems and are not novel by themselves. The main contribution of this paper is the demonstration of their effective use in solving PDE problems in high-dimensional setting.

We now define in detail the architecture we have used. The input $\mathbf x \in \Omega \subset \mathbb R^d$ is the point in the domain at which we wish to calculate the solution while the output $n^{{\rm LSTM}}_\theta(\mathbf x) \in \mathbb R$ is the unnormalized solution at $\mathbf x$. The functions $\left\{\mathtt f_i, \mathtt g_i, \mathtt r_i, \mathtt s_i, \mathtt c_i, \mathtt h_i\right\}$ defined below map $\mathbb R^d$ to $\mathbb R^m$ for $i\in\{1,2,\cdots, L\}$ and they define the hidden layers. The function $\mathtt d_L: \mathbb R^m \to \mathbb R$ is the output layer. For uniformity of notation, we define $\mathtt c_0(\mathbf x)\stackrel{\rm def}{=}\,\mathbf 0_m$ and $\mathtt h_0(\mathbf x)\stackrel{\rm def}{=}\,\mathbf 0_d$ to be the zero vectors of dimension $m$ and $d$ respectively and $\odot$ denotes the Hadamard product.
%
%
%
%
\begin{align}
    \mathtt f_i(\mathbf x) \stackrel{\rm def}{=}& \mathtt A(\mathtt W_f^{(i)}\mathbf x + \mathtt U_f^{(i)}\mathtt h_{i-1}(\mathbf x) + \mathtt b_f^{(i)})\,,\label{eq:layer-f}\\
    \mathtt g_i(\mathbf x) \stackrel{\rm def}{=}& \mathtt A(\mathtt W_g^{(i)}\mathbf x + \mathtt U_g^{(i)}\mathtt h_{i-1}(\mathbf x) + \mathtt b_g^{(i)})\,,\label{eq:layer-g}\\
    \mathtt r_i(\mathbf x) \stackrel{\rm def}{=}& \mathtt A(\mathtt W_r^{(i)}\mathbf x + \mathtt U_r^{(i)}\mathtt h_{i-1}(\mathbf x) + \mathtt b_r^{(i)})\,,\label{eq:layer-r}\\
    \mathtt s_i(\mathbf x) \stackrel{\rm def}{=}& \mathtt A(\mathtt W_s^{(i)}\mathbf x + \mathtt U_s^{(i)}\mathtt h_{i-1}(\mathbf x) + \mathtt b_s^{(i)})\,,\label{eq:layer-s}\\
    \mathtt c_i(\mathbf x) \stackrel{\rm def}{=}&  \mathtt f_i(\mathbf x)\odot \mathtt c_{i-1}(\mathbf x) + \mathtt g_i(\mathbf x)\odot s_i(\mathbf x)\,,\label{eq:layer-c}\\
    \mathtt h_i(\mathbf x) \stackrel{\rm def}{=}& \mathtt r_i(\mathbf x)\odot \mathtt A(\mathtt c_i(\mathbf x))\,,\label{eq:layer-h}\\
    \mathtt d_L(\mathbf y)\stackrel{\rm def}{=}&\mathtt W^\top \mathbf y+ \mathtt b \,, \qquad \mathbf y \in \mathbb R^m\,,\label{eq:layer-final}\\
    n^{{\rm LSTM}}_\theta \stackrel{\rm def}{=}& \mathtt d_L\circ \mathtt h_L\,.\label{eq:def-LSTM} 
\end{align}
The trainable parameters are
\begin{align}
    \theta=\{\mathtt W_f^{(i)}, \mathtt U_f^{(i)}, \mathtt b_f^{(i)}, \mathtt W_g^{(i)}, \mathtt U_g^{(i)}, \mathtt b_g^{(i)}, \mathtt W_r^{(i)}, \mathtt U_r^{(i)}, \mathtt b_r^{(i)}, \mathtt W_s^{(i)}, \mathtt U_s^{(i)}, \mathtt b_s^{(i)}:
i=1,\cdots,L\}\cup\{\mathtt W, \mathtt b\}\,.\label{eq:theta-composition}
\end{align}
The dimensions of these parameters are given below.
\begin{align}
   \mathtt W_f^{(i)}, 
   \mathtt W_g^{(i)},  \mathtt W_r^{(i)},  \mathtt W_s^{(i)} \in&\;
   \mathbb R^{m\times d}\,,\\
   \mathtt U_f^{(i)},
\mathtt U_g^{(i)},
\mathtt U_r^{(i)},
\mathtt U_s^{(i)}\in&
   \begin{cases}\mathbb R^{m\times d},\quad\text{ if }i=1 \\
   \mathbb R^{m\times m},\quad\text{otherwise}
   \end{cases}\\
   \mathtt b_f^{(i)},\mathtt b_g^{(i)},\mathtt b_r^{(i)},\mathtt b_s^{(i)}\in&\;\mathbb R^m\,,\\
   \mathtt W\in\mathbb R^{m}, \mathtt b\in&\;\mathbb R\,.
\end{align}
which implies the size of the network or cardinality of $\theta$ is
\begin{align}
    C=4m[d(L+1)+m(L-1)]+5m+1\label{eq:network-size}\,.
\end{align}
Note that \eqref{eq:network-size} implies that the size of the network grows only linearly with dimension $d$ of the original PDE problem~\eqref{eq:SFPE-0}. This is an important factor for mitigating the curse of dimensionality. We use elementwise $\tanh$ as our activation function,
\begin{align}
    \mathtt A=\tanh\label{eq:activation-choice}\,.
\end{align}
We use $m=50$ and $L=3$ for our experiments which implies our network has $6L=18$ hidden layers. We use the popular Xavier or Glorot initialization \cite{glorot2010understanding}, \cite{datta2020survey} to initialize $\theta$. With that, the description of our architecture is complete. 

\subsection{Optimization}\label{ssec-optimization}
In our experiments we use the ubiquitous Adam optimizer \cite{kingma2014adam} which is often used in the PDE solving literature \cite{han2018solving, zhai2022deep, sirignano2018dgm}. We use a piece-wise linear decaying learning rate. Below $k$ denotes the training iteration and $\delta(k)$ is the learning rate.
\begin{align}
    \delta(k)=\begin{cases}
        5\times10^{-3},\quad\text{if } k<1000
        \,,\\
        1\times10^{-3},\quad\text{if }1000\le k<2000\,,\\
        5\times10^{-4},\quad\text{if }2000\le k<10000\,,\\
        1\times10^{-4},\quad\text{if }k\ge10000\,.\\
        \end{cases}\label{eq:learning-rate}
\end{align}
We stop training after reaching a certain number of iterations $E$ which varies depending on the problem. In all our experiments we use $N=1000$ as the sample size and $\tau=10$ as the resampling interval for algorithm~\ref{algo:steady}.

%% file: steady-tex/results.tex
We are now ready to describe the results of our experiments. Next few sections parallel the examples in section~\ref{sec-examples} and contain problem-specific details about algorithm~\ref{algo:steady} e.g. $\Omega_I, E$ etc. All computations were done with float32 numbers. We comment about the comparison between float32 and float64 computations in the appendix~\ref{ssec-float}.
\subsection{2D ring system}
\begin{figure}[!ht]
    \centering
\includegraphics[scale=0.6
]{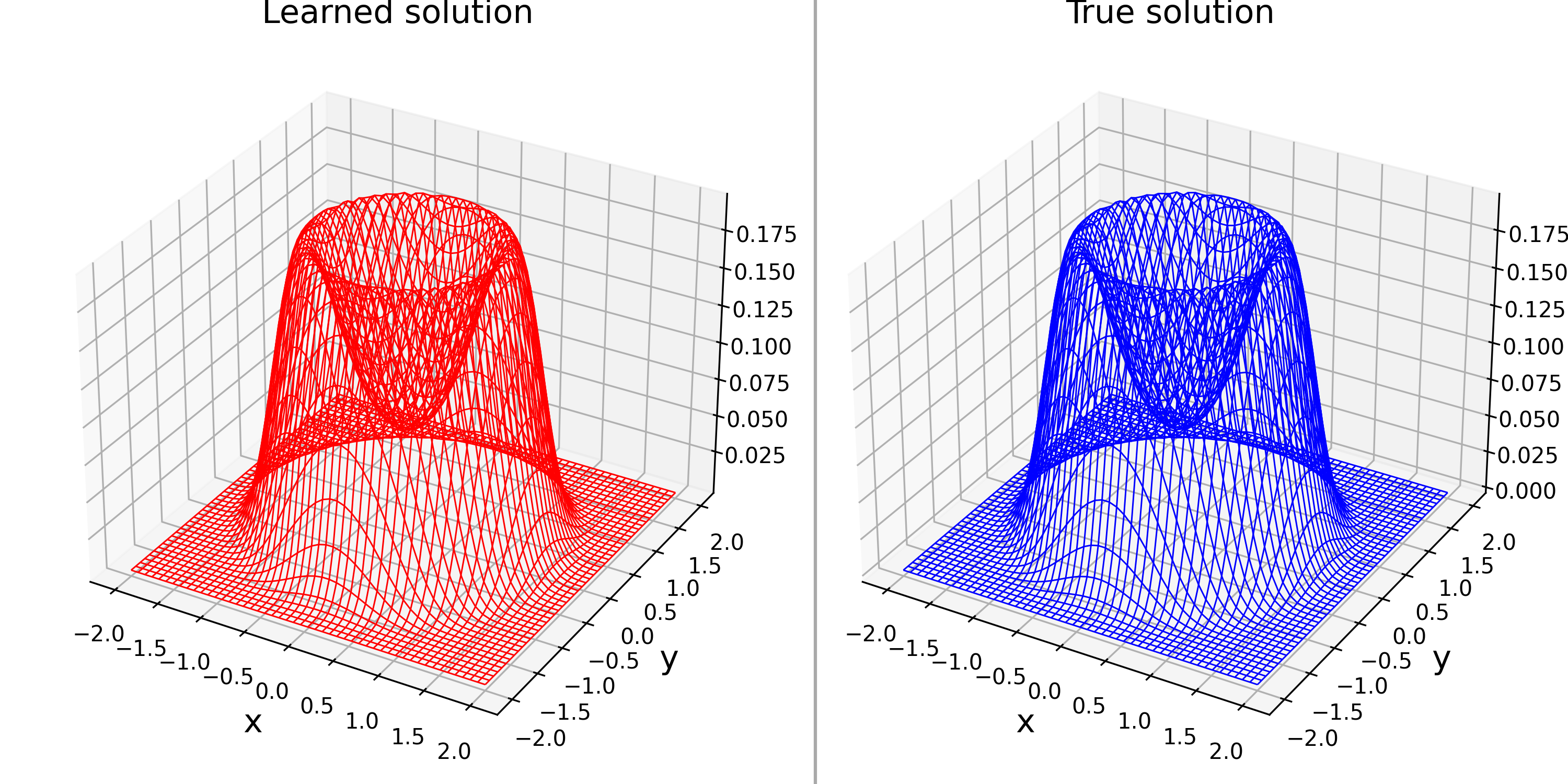}
    \caption{Solution for the 2D ring system}
    \label{fig:2D-surface}
\end{figure}

Figure~\ref{fig:2D-surface} shows the learned and true solutions for the 2D ring system. Note that algorithm~\ref{algo:steady} produces an unnormalized zero of $\mathcal L$ but on the left panel the learned solution has been normalized for easier visualization. 
In this case we use $\Omega_I=[-2,2]^2$ and $E=8\times10^5$ iterations. 
\subsubsection{Comparison with Monte Carlo}\label{sssec-MC-comparison} 
\begin{figure}[!ht]
    \centering
\includegraphics[scale=0.32]{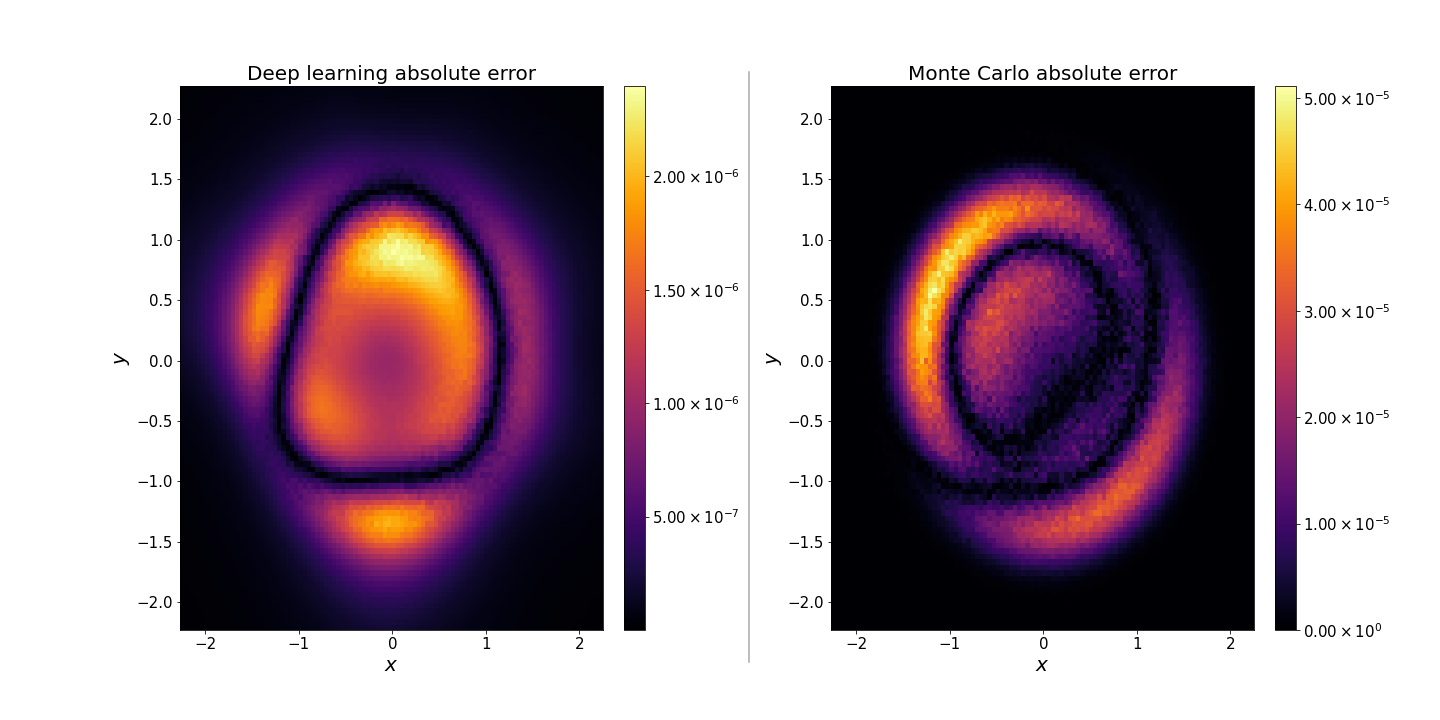}
\caption{Comparison of absolute errors for deep learning and Monte Carlo solutions for the 2D ring system}
    \label{fig:MC-comparison}
\end{figure}

Since the network was trained with domain resampling every $10$ steps and a mini-batch size of $N=1000$, during the entire training procedure $8\times10^7$ points were sampled from the domain. We compute the steady state with Monte Carlo with $8\times10^7$ particles to compare errors produced by both methods. Here the SDE trajectories were generated till time 10 with time-steps of 0.01. Since in this case we know the analytic solution we can compute and compare absolute errors. As we can see in figure~\ref{fig:MC-comparison}, for the same number of overall sampled points, Monte Carlo error is an order of magnitude larger than deep learning error. 

\subsection{2nD ring system}
\begin{figure}[!ht]
    \centering\includegraphics[scale=0.6
]{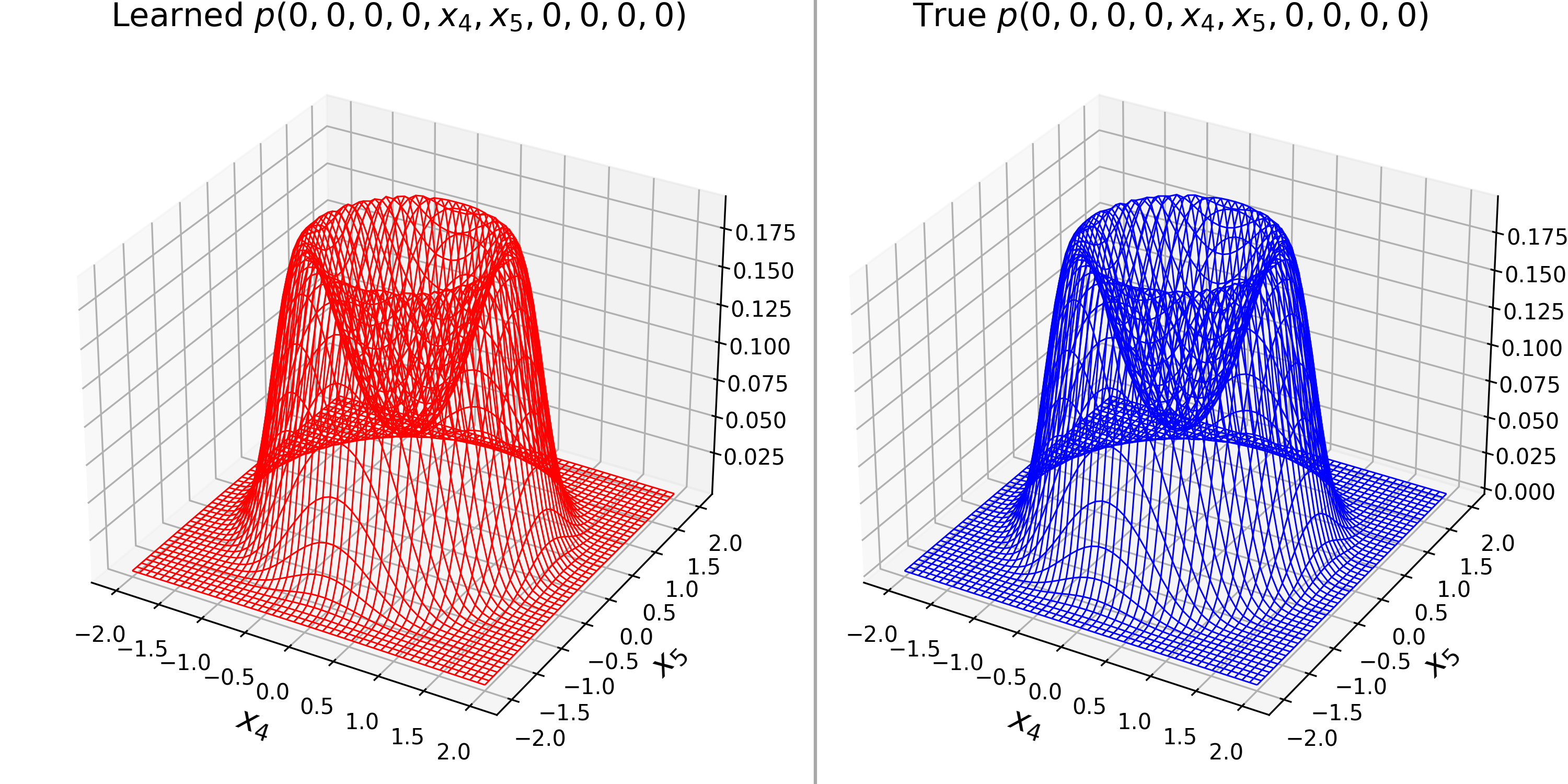}
\caption{Solutions for the 10D ring system. Both solutions have been normalized such that $\int_{\mathbb R}\int_{\mathbb R}p(0, 0, 0, 0, x_4, x_5, 0, 0, 0, 0)\,dx_4\,dx_5=1$} 
    \label{fig:10D-surface}
\end{figure}

\begin{figure}[!ht]
    \centering
\includegraphics[scale=0.32]{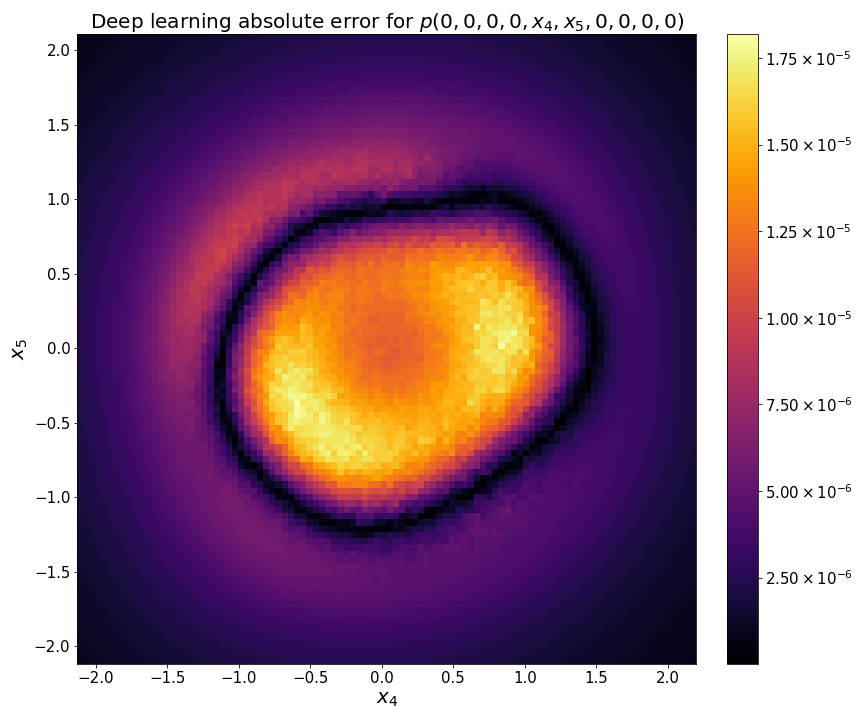}
\caption{Absolute error in the learned solution for the 10D ring system}
    \label{fig:10D-error}
\end{figure}

Although we solve this system for $n=1,2,3,4,5$, in this section we only produce the results for $n=5$ or $d=10$ to avoid repetition. Figure~\ref{fig:10D-surface} shows the solutions for  the 10D ring system for $\Omega_I=[-2, 2]^{10}$ and $E=4.6\times10^6$. In order to visualize the solution we focus on the quantity
$
    p(0, 0, 0, 0, x_4, x_5, 0, 0, 0, 0)
$.
For a visual comparison with the true solution normalization is desirable. But rather than trying to compute a 10-dimensional integral which is a non-trivial problem in itself we can normalize $
    p(0, 0, 0, 0, x_4, x_5, 0, 0, 0, 0)
$ which is much easier to do and due to the decoupled nature of this problem we can expect an identical result as in figure~\ref{fig:2D-surface} which is what we see in figure~\ref{fig:10D-surface}. In both of the panels the solutions have been normalized in a way such that,
$$\int_{\mathbb R}\int_{\mathbb R}p(0, 0, 0, 0, x_4, x_5, 0, 0, 0, 0)\,dx_4\,dx_5=1\,.$$
The error in the learned solution can be seen in figure~\ref{fig:10D-error}.

\subsection{Noisy Lorenz-63 system} 
\begin{figure}[!ht]
    \centering\includegraphics[scale=0.21]{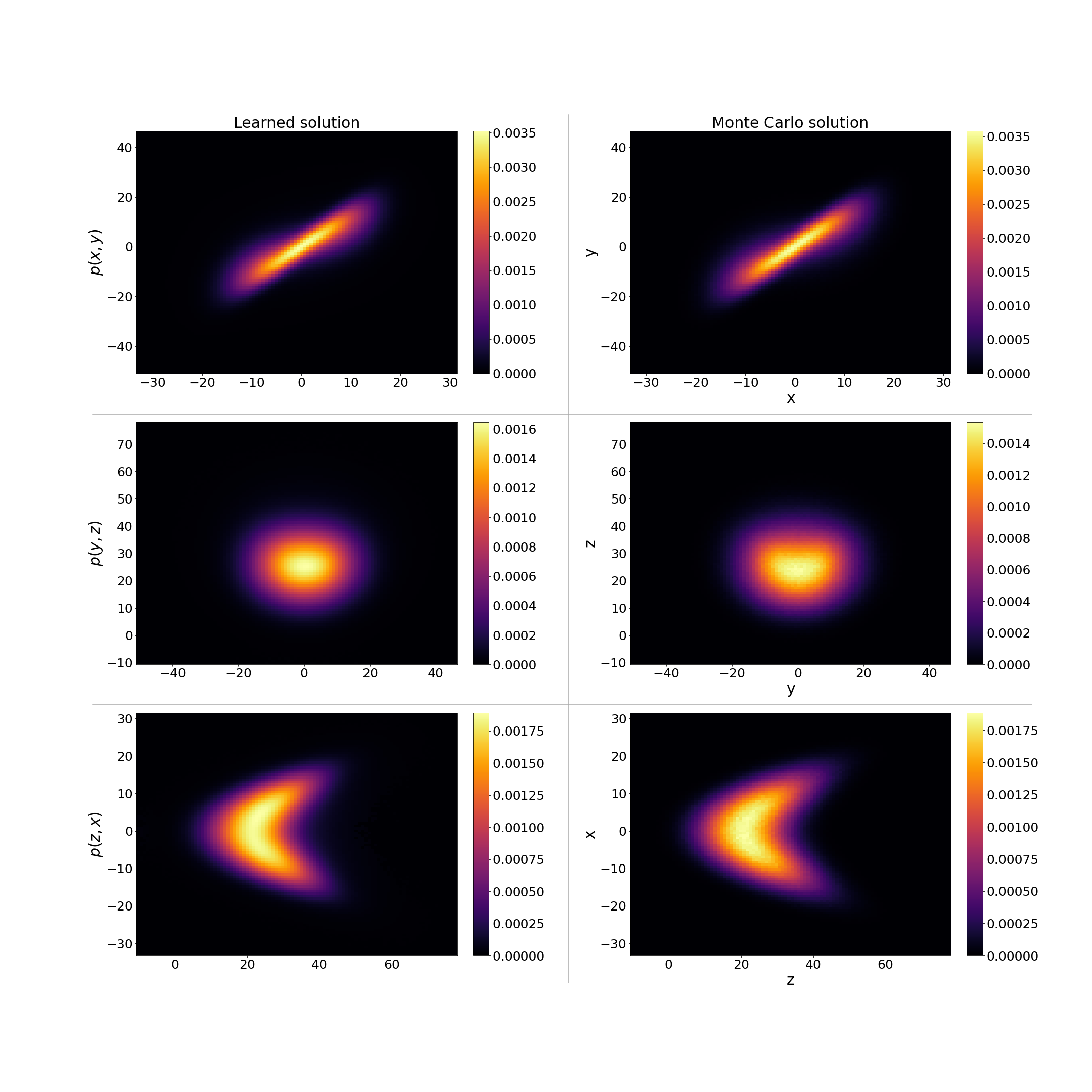}  \caption{Solutions for the noisy Lorenz-63 system}
    \label{fig:L63-steady}
\end{figure}

Figure~\ref{fig:L63-steady} shows the results for the L63 system for $\Omega_I = [-30, 30]\times[-40, 40]\times[0, 70]$ and $E=10^6$. For ease of visualization the solutions have been normalized and in each row one of the dimensions has been integrated over the relevant interval to produce 2D marginals. In order to integrate out one dimension we use a composite Gauss-Legendre quadrature rule. We subdivide the relevant interval into 240 subintervals and use 10-point Gauss-Legendre rule to compute the integral over every subinterval. Note that since $n_\theta$ is a smooth function, our integrand is always a smooth function. The largest possible subinterval is of length $\frac{40-(-40)}{240}=\frac{1}{3}$ so assuming absolute value of the $20$-th derivative of the integrand is upper-bounded by $M$ everywhere, the integration error on each subinterval is upper-bounded by $\frac{2M}{20!}\left(\frac{1}{6}\right)^{20}\le2.25M\times10^{-34}$, see appendix~\ref{ssec-error-GL} for more details on this estimate. To produce the Monte Carlo solution, SDE trajectories were generated till time 10 with time-steps of $10^{-2}$. 
Since Monte Carlo produces lower-accuracy solutions even in lower dimensions as we saw in section~\ref{sssec-MC-comparison} and an analytic solution is unavailable in this case, we cannot produce ``error'' plots with respect true solution, as in the previous example.

\subsection{Noisy Thomas system}
\begin{figure}[!htp]
    \centering\includegraphics[scale=0.32]{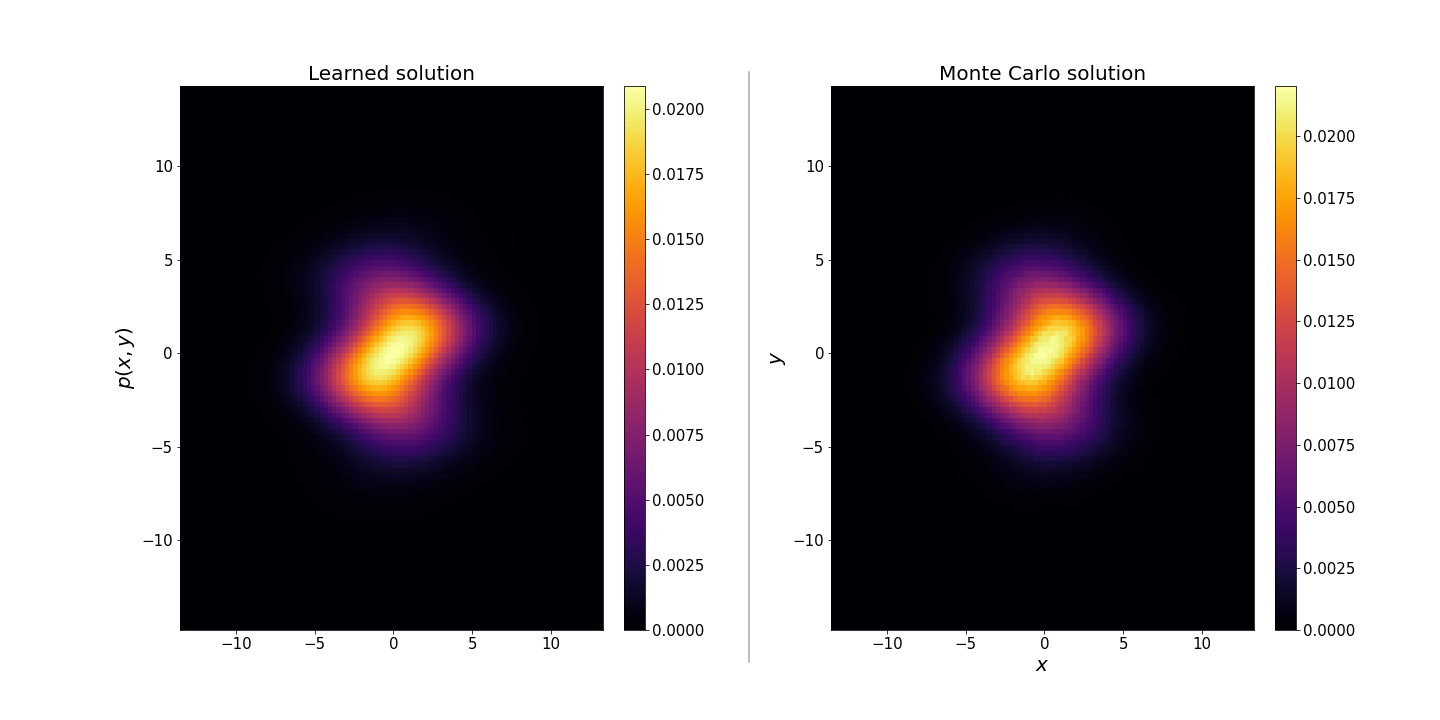}
    \caption{Solutions for the noisy Thomas system}
    \label{fig:Thomas-steady}
\end{figure}

Figure~\ref{fig:Thomas-steady} shows the results for the Thomas system for $\Omega_I = [-10, 10]^3$ and $E=4\times10^5$. Due to the inherent symmetry of this problem it suffices to compute only the 2D marginal $p(x, y)$. To integrate out the $z$ dimension we use 8-point composite Gauss-Legendre quadrature rule with $165$ subintervals. Assuming absolute value of the $16$-th derivative of the integrand is upper-bounded by $M$ everywhere, the integration error on each subinterval is upper-bounded by $\frac{2M}{16!}\left(\frac{10}{165}\right)^{16}\le3.17M\times10^{-33}$, see appendix~\ref{ssec-error-GL} for more details on this error estimate. To produce the Monte Carlo solution, SDE trajectories were generated till time 10 with time-steps of $10^{-2}$. , Thomas system turns out to be the \textit{easiest} among the problems we have solved, even easier than a lower (two) dimensional problem, i.e. algorithm~\ref{algo:steady} converges faster for this system compared to the other ones, as we will see in the next section and in figure~\ref{fig:time-loss}. This could be due to the high degree of symmetry, namely, invariance under permutation of variables and globally Lipschitz drift $\mu$ which is not the case for other systems.
\subsection{Dimension dependence} 
\begin{figure}[!ht]
    \centering\includegraphics[scale=0.4]{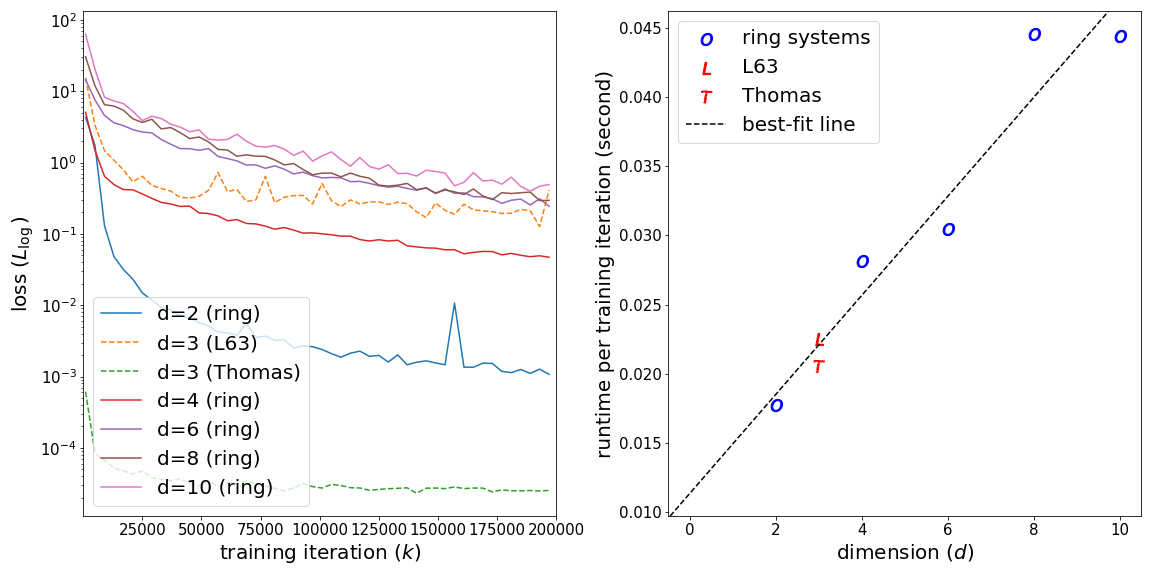}
    \caption{Left panel: Loss vs training iteration starting from iteration $100$. Right panel: time taken per training iteration vs dimension.}
    \label{fig:time-loss}
\end{figure}

In this section we explore the dimension dependence of algorithm~\ref{algo:steady}.
In the left panel of figure~\ref{fig:time-loss} we have plotted the loss given by \eqref{eq:def-steady-loss} against training iterations for all the of the systems above in a semi-log manner starting from iteration $100$. We often encounter spikes in the loss curve for the following reasons
\begin{itemize}
    \item the loss curves are single realizations of algorithm~\ref{algo:steady} instead of being an average
    \item we resample the domain every $10$ iterations and if the new points belong to a previously unexplored region in $\Omega_I$, the loss might increase.
\end{itemize}
But the general trend of loss diminishing with iterations is true for every system. We also see that loss is system-dependent and the \textit{hardness} of these problems or how quickly algorithm~\ref{algo:steady} converges depends on the nature of $\mu$ as much as the dimension. This is easily seen by noting that the two 3D systems (L63 and Thomas) sandwich the 2D and the 4D ring systems in the left panel of figure~\ref{fig:time-loss}. The loss for Thomas system drops very quickly compared to the rest of the systems due to the symmetry and global Lipschitzness of the corresponding drift function. We also see from the right panel of figure~\ref{fig:time-loss} that time taken per training iteration grows near-linearly with dimension. 

We note a couple of points related to the overall time taken for training.
Firstly, since it is hard to estimate the number of iterations required for the loss to drop below a pre-determined level, we refrain from plotting the total runtime of algorithm~\ref{algo:steady} against dimension. In fact, our choice for the number of iterations $E$ for different systems was somewhat ad hoc, based on inspecting the value of the loss function and the solution obtained, and varied from $4\times10^5$ to $4.6\times10^6$ for the lower to higher dimensional problems. It is interesting to note that it is sufficient to increase the number of iterations approximately linearly, instead of needing an exponential growth of this number, with dimension. But a more detailed study and understanding of this aspect as well as the dependence of the required number of training steps on the nature of the drift $\mu$ certainly needs further investigations. 

Second, since the data shown in the right panel of figure~\ref{fig:time-loss} is very much hardware dependent, at this point we note that all of the experiments were done using the cloud service provided by Google Colab. This service automatically assigned runtimes to different hardware depending on availability at the time of computation which explains why the 8D and 10D ring systems take nearly the same amount of time per iteration in figure~\ref{fig:time-loss}.
\subsection{Comparison of loss and distance from truth} 
\begin{figure}[!ht]
    \centering\includegraphics[scale=0.4]{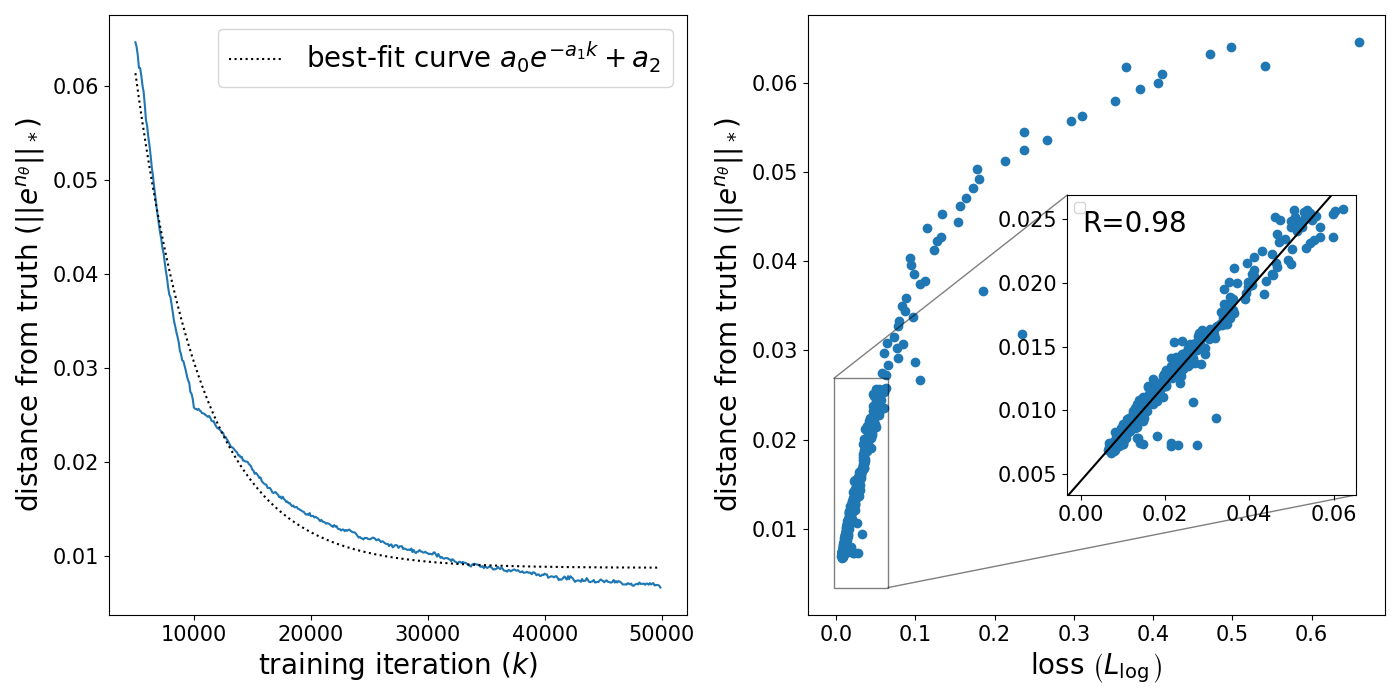}
    \caption{Left panel: Distance from truth vs training iteration every 100 iterations, starting from iteration 5000 and ending at iteration 50000 for the 2D ring system. Right panel: Scatter plot for loss vs distance from truth for the 2D ring system. The inset shows that asymptotically loss and the distance from the truth are linearly related. The inset depicts the data from training iteration $10000$ to $50000$.}
    \label{fig:dist-loss}
\end{figure}

In this section we explore the relationship between the loss given by \eqref{eq:def-steady-loss} and the distance from truth. In spite of being structurally completely different, both are measures of goodness for a computed solution. In most cases we only have access to the loss and therefore it is an important question if a decreasing loss implies getting closer to the truth for algorithm~\ref{algo:steady}. We define the distance of the learned zero from the true solution as follows,
\begin{align}
    \|\phi\|_*\stackrel{\rm def}{=}\sup_{\mathbf x\in\Omega_I} |c\phi(\mathbf x) - p^{\rm true}(\mathbf x)|\label{eq:def-sup-norm}, \qquad c\int_{\mathbb R^d}\phi=1\,,
\end{align}
where $p^{\rm true}$ is the true solution to \eqref{eq:SFPE-0}. \eqref{eq:def-sup-norm} is not easy to compute in arbitrary dimensions but can be computed for the 2D ring system without too much effort since $p^{\rm true}$ is known and the problem is low-dimensional. Figure~\ref{fig:dist-loss} shows the results for the 2D ring system. The right panel of figure~\ref{fig:dist-loss} shows that loss and distance from truth are strongly correlated for algorithm~\ref{algo:steady}. Moreover, asymptotically for small values of the loss function they are linearly related with a Pearson correlation coefficient $R=0.98$ as can be seen from the inset in the right panel which depicts the data from training iteration 10000 to 50000. The best-fit line is also shown in the inset. On the left panel we see that the distance from truth monotonically decreases with training iteration and is extremely well approximated by a curve of the form $a_0e^{-a_1k}+a_2$. Both panels contain data from training iteration 5000 to 50000. We omit the first few iterations to filter out the effects of the random initialization of the trainable parameters. Figure~\ref{fig:dist-loss} serves as a good justification for algorithm~\ref{algo:steady} since it shows that minimizing the loss is akin to getting closer to a true non-trivial zero of $\mathcal L$.
 
%% file: steady-tex/conclusions.tex
In this work we demonstrate the use of deep learning algorithms for finding stationary solutions of the Fokker-Planck equation~\eqref{eq:SFPE-0}. In particular, we find the non-trivial zeros of Fokker-Planck operator $\mathcal L$ defined in~\eqref{eq:SFPE-1}, in the case when the corresponding drift $\mu$ is non-solenoidal. The main motivation is to solve high dimensional Fokker-Planck equations, including the time dependent ones as demonstrated in a sequel~\cite{mandal2024solving}.

We illustrate the method on a variety of problems up to 10 dimensions, with networks whose size (number of parameters) and hence computational time per training iteration both scale linearly with dimension. In all the examples we studied, we notice that it is sufficient for the number of training iterations to grow approximately linearly with dimension in order to obtain similar convergence towards a zero of the cost function, thus leading to overall computational costs that scale roughly quadratically with dimension, which is one of the main advantages of this deep network method compared to other methods.

The rate of convergence of the cost function~\eqref{eq:def-steady-loss} towards zero during training depends on the dimension but also varies significantly with the nature of the problem - some high dimensional problems converge faster than some other low dimensional ones. A more detailed study of this aspect will be an interesting future investigation.

In high dimensions, Monte Carlo methods and the deep network methods such as the one used in this paper are the only viable alternatives, since computational and memory costs of grid based methods scale exponentially making them infeasible. Hence we compare our results with those obtained from Monte Carlo. In examples where analytical solutions are known, it is seen that the deep network solutions are more accurate than Monte Carlo solutions obtained with similar computational cost. The other main advantage of the deep network method is that we can get solutions in a functional form which the Monte Carlo is incapable of doing.

The deep network is trained by minimizing the loss~\eqref{eq:def-steady-loss} which by itself does not imply {\it a priori} that the solution obtained gets closer to the true desired solution. For problems for which analytical solutions are available, we look at the relation between the loss and the distance of the function represented by the network from the true (analytical) zero of the Fokker-Planck operator. Even though these quantities are structurally completely different, we see that they are strongly correlated. Moreover, they can be asymptotically linearly related for small values of the loss function.

The results in this work lead to several possible avenues for further investigations. We have already mentioned some of these, such as a deeper understanding of the relation between the loss function and the distance from the true solution as well as the number of training iterations required for obtaining a pre-determined level of the loss. An additional interesting direction is to explore the geometric questions related to the landscape of the loss defined in~\eqref{algo:steady}. For example, in case the nullspace of $\mathcal L$ is 1-dimensional, as is the case for the problems considered in this paper, it will be a challenging problem to understand the relation of this nullspace to the set of all the minima $\theta^*$ of $L_{\log}$ defined in~\eqref{eq:def-steady-loss}, including the topological properties of these sets.

%% file: steady-tex/appendix.tex
\subsection{Existence and uniqueness of solutions to example problems}\label{ssec-unique} In this section we prove that the example problems used here have a unique weak solution in $W^{1,2}_{\rm loc}(\mathbb R^d)$.
 We employ the method of Lyapunov function as described in \cite{huang2015steady} to arrive at existence and uniqueness. First we begin with the prerequisites for this approach.
 \subsubsection{Lyapunov functions}
\begin{defn}
    Let $U \in C(\mathcal U)$ be a non-negative function and denote $\rho_M = \sup_{\mathbf x\in \mathcal U} U(\mathbf x)$,
called the essential upper bound of $U$. $U$ is said to be a compact function in $U$ if
\begin{align}
i)\; U (\mathbf x) < \rho_M,\quad \mathbf x \in \mathcal U
\end{align}
and
\begin{align}
ii) \lim_{\mathbf x\to\partial U} U (\mathbf x) = \rho_M\,.
\end{align}
\end{defn} This definition of a compact function appears as definition 2.2 in \cite{huang2015steady}.
\begin{prop}
An unbounded, non-negative function $U\in C(\mathbb R^d)$ is compact iff 
\begin{align}
    \lim_{\|\mathbf \mathbf x\|_2\to+\infty} U(\mathbf x) = +\infty\,.
\end{align}
\end{prop}
This proposition appears as proposition 2.1 in \cite{huang2015steady}.

\begin{defn}
    Let $U$ be a compact function in $C^2(\mathcal U)$ with essential upper bound $\rho_M$. $U$ is called a Lyapunov function in $\mathcal U$ with respect to $\mathcal L^*$ is $\exists\,\rho_m\in(0, \rho_M)$ and a constant $\gamma>0$ such that
    \begin{align}
        \mathcal L^* U(\mathbf x)\le -\gamma,\qquad\forall \mathbf x\in \mathcal U\setminus\overline{\{\mathbf x\in\mathcal U: U(\mathbf x)<\rho_m\}}\,,
    \end{align}
    where $\mathcal L^*$ is the adjoint Fokker-Planck operator given by
    \begin{align}
        \mathcal L^*f = \mu\cdot\nabla f+ D\odot\nabla^2 f\,. \label{eq:def-adjoint-FP-op}
    \end{align}
\end{defn}
This definition appears as definition 2.4 in \cite{huang2015steady}.
Now we are ready to state the main theorem that will help us prove uniqueness for our example problems.
\begin{thm}
    If the components of $\mu$ are in $L^2_{\rm loc}(\mathcal U)$ and there exists a Lyapunov function with respect to $\mathcal L^*$ in $C^2(\mathcal U)$ then \eqref{eq:SFPE-0} has a positive weak solution in the space $W^{1, 2}_{\rm loc}(\mathcal U)$. If, in addition,
the Lyapunov function is unbounded, the solution is unique in $\mathcal U$.\label{thm:unique-steady}
\end{thm}
This theorem appears as theorem $A$ in \cite{huang2015steady}. Since the components of $\mu$ are locally integrable for our example problems, all we need to do is find an unbounded Lyapunov function $U$ for proving existence and uniqueness in $W^{1,2}_{\rm loc}(\mathbb R^d)$.
\subsubsection{Existence and uniqueness of solution for 2D ring system}\label{sssec-2D-unique}
Setting 
\begin{align}
\mathcal U &= \mathbb R^2\,,\\
U(x, y) &= x^2+y^2\,,\\
\rho_m &= \frac{1}{2}+\sqrt{D+1}\,,\\
\gamma &= 4D+6\,,\\
\end{align}
we see that,
\begin{align}
    \mathcal L^*U +\gamma = -8\left(x^2+y^2-\frac{1}{2}\right)^2 + 8(D+1)
\end{align}
and
\begin{align}
    \mathcal U\setminus\overline{\{\mathbf x\in\mathcal U: U(\mathbf x)<\rho_m\}} = \left\{(x,y)\in\mathbb R: x^2+y^2>\rho_m\right\}\,.
\end{align}
In  $\left\{(x,y)\in\mathbb R: x^2+y^2>\rho_m\right\}$, 
\begin{align}
    \mathcal L^*U+\gamma \le 0
\end{align}
and therefore $U$ is an unbounded Lyapunov function for the 2D ring system which guarantees uniqueness of solution \eqref{eq:grad-sol}.

\subsubsection{Existence and uniqueness of solution for L63 system}\label{sssec-L63-unique}
Setting,
\begin{align}
U(x, y, z) = \rho x^2 +\alpha y^2 + \alpha(z-2\rho)^2\,,
\end{align}
we see that
\begin{align}
    \mathcal L^*U &= -2\alpha\rho x^2 - 2\alpha y^2 -2\alpha\beta z^2 + 4\alpha\beta\rho z + 2D(2\alpha+\rho)\,,\\
    &=-2\alpha\rho x^2 - 2\alpha y^2 -\alpha\beta z^2 -\alpha\beta(z-2\rho)^2 + 4\alpha\beta\rho^2 + 2D(2\alpha+\rho)\,,\\
    &\le -\rho x^2 -\alpha y^2 -\alpha(z-2\rho)^2 + 4\alpha\beta\rho^2 + 2D(2\alpha+\rho)\,,\label{eq:L63-params-bigger-than-1}\\
    &= -U(x, y, z)+ 4\alpha\beta\rho^2 + 2D(2\alpha+\rho)\,.
\end{align}
\eqref{eq:L63-params-bigger-than-1} is a consequence of $\alpha, \beta, \rho>1$. Now setting,
\begin{align}
    \gamma &= 1,\\
    \rho_m &= 4\alpha\beta\rho^2 + 2D(2\alpha+\rho)+1\,,
\end{align} we see that in $\{U>\rho_m\}$,
\begin{align}
    \mathcal L^*U +\gamma \le 0\,.
\end{align}
So $U$ is an unbounded Lyapunov function for this system and we have a unique solution.
\subsubsection{Existence and uniqueness of solution for Thomas system}\label{sssec-Thomas-unique}
Setting,
\begin{align}
    U(x, y, z) = x^2+y^2+z^2\,,
\end{align}
we see that
\begin{align}
    \mathcal L^* U &= x\sin y + y\sin z + z\sin x - b(x^2+y^2+z^2) + 6D\,,\\
    &\le \sqrt{3U}-bU + 6D\,,\label{eq:CS-on-Thomas-unique}\\
    &= -b\left(\sqrt{U}-\frac{\sqrt{3}}{2b}\right)^2 +\frac{3}{4b}+6D\,.
\end{align}
\eqref{eq:CS-on-Thomas-unique} follows from Cauchy Schwarz inequality. Setting,
\begin{align}
    \gamma &= \frac{1}{4b},\\
    \rho_m &= \left(\frac{\sqrt{3}}{2b}+\frac{\sqrt{1+6bD}}{b}\right)^2\,,
\end{align}
we see that in $\{U>\rho_m\}$,
\begin{align}
    \mathcal L^*U +\gamma \le 0\,.
\end{align}
So $U$ is an unbounded Lyapunov function for this system and we have a unique solution.

\subsection{Monte Carlo steady state algorithm}\label{ssec-MC-algo}
The time-dependent FPE given by
\begin{equation}
\begin{aligned}
    &\frac{\partial  p(t, \mathbf x)}{\partial t} =\mathcal L p(t, \mathbf x),\qquad\mathbf x\in\mathbb R^d,\; t\ge0\,,\\&p(0, \mathbf x)=p_0(\mathbf x),\qquad\mathbf x\in\mathbb R^d\,,\\
    &\int_{\mathbb R^d}p(t,\mathbf x)\,d\mathbf x = 1,\qquad\forall\;t\ge0\,,
    \label{eq:FPE-0}
\end{aligned}
\end{equation}
gives us the probability density of the random process $X_t$ which is governed by the SDE,
\begin{equation}
\begin{aligned}
    &dX_t=\mu\,dt+\sigma\,dW_t\,,\\
    &X_0\sim p_0\,.\label{eq:SDE-0}
\end{aligned}
\end{equation}where $\{W_t\}$ is the standard Wiener process, see for example chapters 4, 5 of \cite{gardiner2009stochastic}. We can evolve $\eqref{algo:steady}$ up to sufficiently long time using Euler-Maruyama method \cite{kloeden1992stochastic} to approximate the steady state solution of \eqref{eq:FPE-0} or the solution of \eqref{eq:SFPE-0} as follows. Here $\mathcal N$ denotes the multivariate normal distribution.
\begin{algorithm}[!ht]
Sample $\{ X_0^{(i)}\}_{i=1}^N\sim p_0$.\\
Set the time-step $h$.\\
Set the number of steps $S$.\\
\For {$k=1,2\cdots, S$}{
 Sample $w^i_k\sim\mathcal N(\mathbf 0_d, h I_d)\;\;\forall\;i$\\
 $ X_k^{(i)}\leftarrow  X_{k-1}^{(i)} + \mu\left(X_{k-1}^{(i)}\right)h + \sigma w_k^i\;\;\forall\;i$\\
}
Subdivide the domain of interest $\Omega_I$ into $d$-dimensional boxes.\\ Count the number of $X^{(i)}_{S}$ that are in a box to estimate the stationary density at the center of the box.
\caption{Monte Carlo steady state algorithm}\label{algo:MC}
\end{algorithm}
Note that in case of a unique solution of \eqref{eq:SFPE-0}, many choices of $p_0$ can lead to the stationary solution. In all our examples, it suffices to choose $p_0$ to be the standard $d$-dimensional normal distribution.

\subsection{System dependence of floating point errors}\label{ssec-float}
As mentioned in section~\ref{ssec-rationale}, floating points errors are an important aspect in scientific computations. In this section we discuss why using higher precision floats may be necessary while solving Fokker-Planck equations with alogrithm~\ref{algo:steady}. We do so by presenting results for two gradient systems, one of which necessitates the use of 64-bit precision while the other one does not, in a way that is made precise below.

\begin{figure}
    \centering
\includegraphics[scale=0.65]{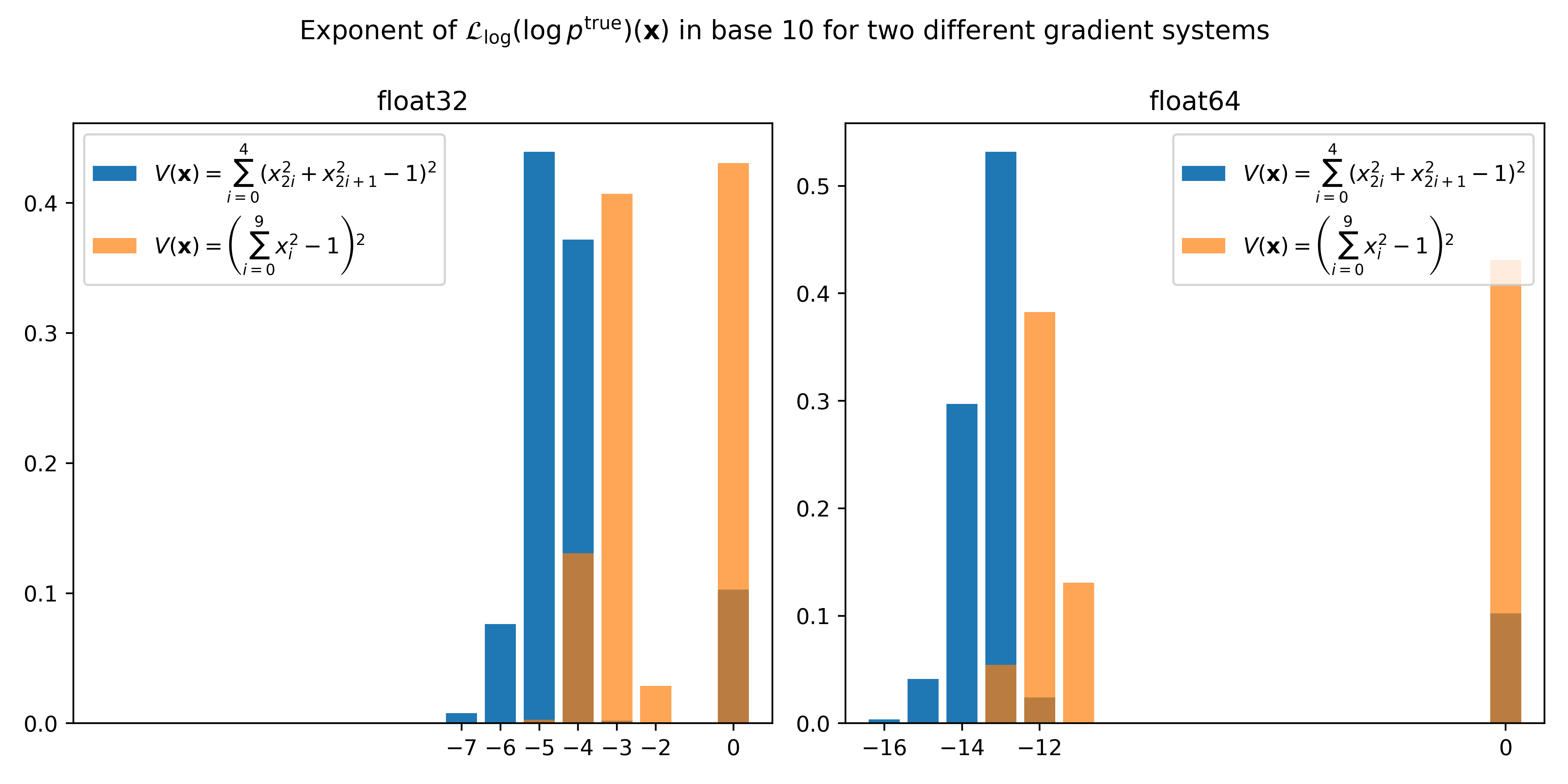}
    \caption{Normalized histogram of the normalized base-10 exponent $b_j$ of $\mathcal L_{\log}(\log p^{\rm true})(\mathbf x_j)$ as in~\eqref{eq:mantissa-exponent}, where $\{\mathbf x_j\}_{j=1}^{10^6}$ is a uniform sample from $[-2, 2]^{10}$. In the left and right panels the computations were done using float32 and float64 numbers respectively. All $4$ histograms have a local maximum around $0$ which corresponds to the samples where $\mathcal L_{\log}(\log p^{\rm true})(\mathbf x_j)$ was numerically evaluated to 0 up to machine precision.} 
    \label{fig:float-histogram}
\end{figure}

The first system is the 10D ring system described in section~\ref{ssec-2nD-ring} with the potential~\eqref{eq:2nD-V}.
The second system has a 10D hypersphere attractor and is given by the following potential, 
\begin{align}
    V(\mathbf x)=\left(\sum_{i=0}^{9}x_{i}^2-1\right)^2\,.\label{eq:near-zero}
\end{align}
We set $D=1$ or $\sigma=\sqrt 2$ for both systems. For each of these systems, since we know the true solution $p^{\rm true}$, we can calculate the following quantity with automatic differentiation: 
\begin{align}
    \mathcal L_{\log}(\log p^{\rm true})(\mathbf x_j)=\mathcal L_{\log}\left(-\frac{V}{D}\right)(\mathbf x_j),\qquad j=1,2,\cdots N\,,
\end{align}
where $\{\mathbf x_j\}_{j=1}^N$ is a uniform sample from $[-2, 2]^{10}$. Analytically the quantity appearing in~\eqref{eq:near-zero} is $0$ and when evaluated numerically, we expect either zero or near zero numbers. Suppose, when written in normalized mantissa-exponent form in base $10$~\cite{schmid1974decimal}, this quantity looks like
\begin{align}
    \mathcal L_{\log}(\log p^{\rm true})(\mathbf x_j) =  a_j\times10^{b_j}\,.\label{eq:mantissa-exponent}
\end{align}
We adopt the convention that $b_j=0$ when the LHS of~\eqref{eq:mantissa-exponent} is $0$. Therefore, we would expect $b_j$ to be either $0$ when the associated float is $0$ up to machine precision or a highly negative integer as allowed by the corresponding floating point system. Figure~\ref{fig:float-histogram} shows the normalized histogram of $b_j$ for both systems for float32 and float64 for sample size $N=10^6$ in each case. Both systems have peaks at $0$ in both float32 and float64 which correspond to the samples where the quantity in~\eqref{eq:near-zero} numerically evaluates to $0$ up to machine precision. This accounts for nearly half the samples for the second system and only about 10\% of the samples for the first system. But for the samples for which $b_j\neq0$, $\mathcal L_{\log}(\log p^{\rm true})(\mathbf x_j)$ evaluates much closer to $0$ for the first system when compared to the second system. This deviation from zero for the second system is more prominent when the computation is done in float32 with the majority of the nonzero samples having $b_j=-3$. For the same floating point system, the samples for the first system that are furthest from $0$ have $b_j=-4$ and the majority of the nonzero samples have $b_j=-5$. This indicates that float64 might be a more appropriate choice for the second system. A more detailed study of this dependence of the results on the choice of the floating point precision for a variety of systems is an interesting avenue for future research.

\subsection{Integration error for \texorpdfstring{$n$}{Lg}-point Gauss-Legendre rule}\label{ssec-error-GL}
Suppose we are trying to integrate a smooth function $f(x)$ over $\left[a-\frac{h}{2}, a+\frac{h}{2}\right]$ with $n$-point Gauss-Legendre rule where $h\in(0, 1]$. Let us denote $I[f]$ to be the Gauss-Legendre approximation of $\int_{a-\frac{h}{2}}^{a+\frac{h}{2}} f(x)\,dx$. Recalling that $n$-point Gauss-Legendre gives us exact integrals for polynomial of degree $\le 2n-1$ and using the Lagrange form of Taylor remainder we see that 
\begin{align}
    \left|I[f]-\int_{a-\frac{h}{2}}^{a+\frac{h}{2}} f(x)\,dx\right| &\le MI\left[\frac{(x-a)^{2n}}{(2n)!}\right]+M\int_{a-\frac{h}{2}}^{a+\frac{h}{2}} \frac{(x-a)^{2n}}{(2n)!}\,dx\,,\label{eq:diff-GL-true}
\end{align}
where $|f^{(2n)}(x)|\le M\;\;\forall\;\;x\in\left[a-\frac{h}{2}, a+\frac{h}{2}\right]$. To bound the first term on the RHS of \eqref{eq:diff-GL-true} we can use the fact that if 
\begin{align}
    I[f] = \sum_{i=1}^nw_if(x_i)\,,
\end{align}
then
\begin{align}
    &I[1] = \int_{a-\frac{h}{2}}^{a+\frac{h}{2}} 1 \,dx = h\,,\\
    \implies&\sum_{i=1}^nw_i = h\le1\,,\\
    \implies&I\left[\frac{(x-a)^{2n}}{(2n)!}\right]\le\frac{1}{(2n)!}\left(\frac{h}{2}\right)^{2n}\,.
\end{align}
Therefore,
\begin{align}
    \left|I[f]-\int_{a-\frac{h}{2}}^{a+\frac{h}{2}} f(x)\,dx\right|\le \frac{M}{(2n)!}\left(\frac{h}{2}\right)^{2n}+\frac{2M}{(2n+1)!}\left(\frac{h}{2}\right)^{2n+1}\le\frac{2M}{(2n)!}\left(\frac{h}{2}\right)^{2n}\,.
\end{align}

%% file: steady-tex/acknowledge.tex
This work was supported by the Department of Atomic Energy, Government of India, under project no. RTI4001.